\documentclass[journal, twoside]{IEEEtran}

\usepackage{import}

\usepackage[utf8]{luainputenc}

\usepackage[T1]{fontenc}

\usepackage{graphicx}
\usepackage[dvipsnames]{xcolor}

\usepackage[english]{babel}
\addto\captionsenglish{}
\addto\captionsenglish{}
\usepackage{csquotes}

\usepackage{newtxtext}
\usepackage{amsthm}
\usepackage[slantedGreek]{newtxmath}
\usepackage[OMLmathsfit]{isomath}
\DeclareMathAlphabet{\mathbfsf}{\encodingdefault}{\sfdefault}{bx}{n}
\usepackage{bm}
\usepackage{envmath}
\usepackage{mathtools}
\usepackage{commath}
\usepackage{siunitx}

\usepackage{subcaption}
\makeatletter
\renewcommand{\p@subfigure}{\thefigure}

\makeatother

\usepackage{booktabs}
\usepackage{footmisc}  

\usepackage{url}

\theoremstyle{definition}

\theoremstyle{plain}

\theoremstyle{remark}

\usepackage{lineno}
\modulolinenumbers[5]
\usepackage{todonotes}
\usepackage{umoline}

\usepackage{pgfplots}
\usepackage{pgfplotstable}
\pgfplotsset{compat=newest}
\pgfplotsset{plot coordinates/math parser=false}
\newlength\figureheight
\newlength\figurewidth
\pgfplotsset{every axis plot/.append style={line width=1.5pt},
    legend style={font=\footnotesize, 
        text height=1.0ex,
        draw=black,
        fill=white,
        legend cell align=left}}

\usepackage[hidelinks]{hyperref} 
\usepackage[english]{cleveref}

\Crefname{defn}{definition}{definitions}
\Crefname{defn}{Definition}{Definitions}

\Crefname{asm}{assumption}{assumptions}
\Crefname{asm}{Assumption}{Assumptions}

\crefname{lem}{lemma}{lemmas} 
\Crefname{lem}{Lemma}{Lemmas}

\crefname{prop}{proposition}{propositions} 
\Crefname{prop}{Proposition}{Propositions}

\crefname{thm}{theorem}{theorms} 
\Crefname{thm}{Theorem}{Theorms}

\crefname{cor}{corollary}{corollaries}
\Crefname{cor}{Corollary}{Corollaries}

\Crefname{figure}{Fig.}{Figs.}

\AtBeginEnvironment{appendices}{\crefalias{section}{appendix}}

\usepackage{upgreek}
\newcounter{subequation}
\newlength\mtabskip\mtabskip=-1.25cm

\def\mtabLong{long}

\newcommand{\mr}{\mathrm}

\newcommand{\mc}{\mathcal}

\newcommand{\veg}[1]{\bm{#1}}     
\newcommand{\mat}[1]{\mathsfbfit{#1}} 
\renewcommand{\vec}[1]{\mathsfbfit{#1}} 
\newcommand{\wvec}[1]{\widetilde{\mathsfbfit{#1}}} 
\newcommand{\op}[1]{\mathcal{#1}} 
\newcommand{\vecop}[1]{\bm{\mathcal{#1}}} 

\newcommand{\dd}{\mathrm{d}}  

\newcommand{\jm}{\mathrm{j}}  

\newcommand{\e}{\mathrm{e}}

\DeclareMathOperator{\Div}{div}

\DeclareMathOperator{\grad}{\mathbf{grad}}

\DeclareMathOperator{\supp}{supp}

\newcommand{\T}{\mr{T}}

\newcommand\restr[2]{{
        \left.\kern-\nulldelimiterspace 
        #1 
        \vphantom{|} 
        \right|_{#2} 
}}

\newcommand\rst[3]{{
        \left.\kern-\nulldelimiterspace 
        #1 
        \vphantom{|} 
        \right|_{#2}^{#3} 
}}

\newcommand{\numSplines}{s}     
 
\usepackage{acro}
\usepackage{pdfcomment}
\acsetup{pdfcomments/use=true}

\DeclareAcronym{DG}
{
    short = DG ,
    long = discontinuous Galerkin
}

\DeclareAcronym{ACA}
{
    short = ACA ,
    long = adaptive cross approximation
}

\DeclareAcronym{EFIE}
{
    short =  EFIE ,
    long = electric field integral equation
}

\DeclareAcronym{MFIE}
{
    short =  MFIE ,
    long = magnetic field integral equation
}

\DeclareAcronym{CFIE}
{
    short =  CFIE ,
    long = combined field integral equation
}

\DeclareAcronym{MUIE}
{
    short =  MUIE ,
    long = Müller integral equation
}

\DeclareAcronym{PMCHWT}
{
    short =  PMCHWT ,
    long = Poggio-Miller-Chang-Harrington-Wu-Tsai integral equation
}

\DeclareAcronym{MLFMA}
{
    short =  MLFMA ,
    long = multilevel fast multipole algorithm
}

\DeclareAcronym{MLFMM}
{
    short =  MLFMM ,
    long = multilevel fast multipole method
}

\DeclareAcronym{AIM}
{
    short =  AIM ,
    long = adaptive integral method
}

\DeclareAcronym{BEM}
{
    short =  BEM ,
    long = boundary element method
}

\DeclareAcronym{ACE}
{
    short =  ACE ,
    long = accelerated Cartesian expansion,
    }

\DeclareAcronym{SPD}
{
    short =  SPD ,
    long = {symmetric, positive definite}
}

\DeclareAcronym{SPSD}
{
    short =  SPD ,
    long = {symmetric, positive semi-definite}
}

\DeclareAcronym{PEC}
{
    short =  PEC ,
    long = perfectly electrically conducting
}

\DeclareAcronym{RWG}
{
    short = RWG ,
    long = Rao-Wilton-Glisson
} 

\DeclareAcronym{BC}
{
    short = BC ,
    long = Buffa-Christiansen
}

\DeclareAcronym{SVD}
{
    short = SVD ,
    long = singular value decomposition
}

\DeclareAcronym{CG}
{
    short = CG ,
    long = conjugate gradient
} 

\DeclareAcronym{PCG}
{
    short = PCG ,
    long = preconditioned conjugate gradient
} 

\DeclareAcronym{CGS}
{
    short = CGS ,
    long = conjugate gradient squared
}

\DeclareAcronym{CMP}
{
    short = CMP ,
    long = Calderón multiplicative preconditioner
} 

\DeclareAcronym{RFCMP}
{
    short = RF-CMP ,
    long = refinement-free Calderón multiplicative preconditioner
} 

\DeclareAcronym{HPD}
{
    short = HPD ,
    long = {Hermitian, positive definite}
} 

\DeclareAcronym{RHS}
{
    short = RHS ,
    long = right-hand side
}

\DeclareAcronym{LSE}
{
    short = LSE ,
    long = linear system of equations ,
    long-plural-form = linear systems of equations 
}

\DeclareAcronym{AMG}
{
    short = AMG ,
    long = algebraic multigrid
}

\DeclareAcronym{PW}
{
    short = PW ,
    long = plane wave
}

\DeclareAcronym{GMRES}
{
    short = GMRES ,
    long = generalized minimum residual
}

\DeclareAcronym{IDR}
{
    short = IDR ,
    long = induced dimension reduction
}

\DeclareAcronym{BICGstab}
{
    short = BiCGstab ,
    long = stabilized bi-conjugate gradient
}

\DeclareAcronym{FF}
{
    short = FF ,
    long = far field
}

\DeclareAcronym{NF}
{
    short = NF ,
    long = near field
}

\DeclareAcronym{SC}
{
    short = SCC ,
    long = standard convergence criterion ,
    pdfcomment = standard convergence criterion
}

\DeclareAcronym{RC}
{
    short = RSCC ,
    long = random-sampling convergence criterion ,
    pdfcomment = random-sampling convergence criterion
}

\DeclareAcronym{CC}
{
    short = CCC ,
    long = combined convergence criterion ,
    pdfcomment = combined convergence criterion
}

\DeclareAcronym{HCA}
{
    short = HCA,
    long = hybrid cross approximation
}

\DeclareAcronym{FMM}
{
    short = FMM ,
    long = fast multipole method  ,
    pdfcomment = fast multipole method 
}

\DeclareAcronym{HMatrix}
{
    short = $\mathscr{H}$-matrix ,
    long = hierarchical matrix ,
    pdfcomment = hierarchical matrix 
}

\DeclareAcronym{MLMDA}
{
    short = MLMDA ,
    long = multilevel matrix decomposition algorithm ,
    pdfcomment = multilevel matrix decomposition algorithm
}

\DeclareAcronym{NCA}
{
    short = NCA ,
    long = nested cross approximation ,
    pdfcomment = nested cross approximation
}

\DeclareAcronym{MVP}
{
    short = MVP,
    long = matrix--vector product ,
    pdfcomment = matrix--vector product 
}

\DeclareAcronym{FFT}
{
    short = FFT,
    long = fast Fourier transform ,
    pdfcomment = fast Fourier transform 
}

\DeclareAcronym{CAD}
{
    short = CAD,
    long = computer aided design,
} 

\newcolumntype {n}{c}
\newcolumntype {N}{>{\small}c}
\newcolumntype {L}{>{\small}l}
\newcolumntype {F}{>{\footnotesize}c}
\newcolumntype {v}[1]{>{\raggedright \hspace {0pt}} p {#1}}
\newcolumntype {V}[1]{>{\small \raggedright \hspace {0pt}} p {#1}}
\newcolumntype{d}[1]{>{\DC@{.}{.}{#1}}c<{\DC@end}}

\newcolumntype{R}[1]{%
    >{\begin{turn}{90}\begin{minipage}{#1}\small\raggedright\hspace{0pt}}l%
            <{\end{minipage}\end{turn}}%
}

\pgfplotsset{colormap={hawaii}{
		rgb = (0.550541, 0.006842, 0.451980);
		rgb = (0.551494, 0.015367, 0.447972);
		rgb = (0.552426, 0.023795, 0.443998);
		rgb = (0.553328, 0.032329, 0.440021);
		rgb = (0.554227, 0.041170, 0.436063);
		rgb = (0.555098, 0.049286, 0.432125);
		rgb = (0.555948, 0.056667, 0.428188);
		rgb = (0.556797, 0.063525, 0.424272);
		rgb = (0.557619, 0.069970, 0.420377);
		rgb = (0.558415, 0.076028, 0.416509);
		rgb = (0.559210, 0.081936, 0.412663);
		rgb = (0.559991, 0.087507, 0.408823);
		rgb = (0.560746, 0.092811, 0.405012);
		rgb = (0.561495, 0.098081, 0.401237);
		rgb = (0.562235, 0.103128, 0.397471);
		rgb = (0.562954, 0.108005, 0.393736);
		rgb = (0.563663, 0.112872, 0.390025);
		rgb = (0.564355, 0.117530, 0.386344);
		rgb = (0.565032, 0.122122, 0.382698);
		rgb = (0.565709, 0.126681, 0.379074);
		rgb = (0.566380, 0.131171, 0.375474);
		rgb = (0.567037, 0.135542, 0.371905);
		rgb = (0.567679, 0.139872, 0.368378);
		rgb = (0.568312, 0.144198, 0.364861);
		rgb = (0.568939, 0.148416, 0.361384);
		rgb = (0.569559, 0.152618, 0.357942);
		rgb = (0.570171, 0.156806, 0.354519);
		rgb = (0.570777, 0.160934, 0.351127);
		rgb = (0.571377, 0.165008, 0.347764);
		rgb = (0.571972, 0.169120, 0.344417);
		rgb = (0.572562, 0.173131, 0.341120);
		rgb = (0.573142, 0.177166, 0.337836);
		rgb = (0.573711, 0.181138, 0.334602);
		rgb = (0.574276, 0.185151, 0.331356);
		rgb = (0.574840, 0.189095, 0.328170);
		rgb = (0.575406, 0.193035, 0.324992);
		rgb = (0.575967, 0.196978, 0.321854);
		rgb = (0.576518, 0.200854, 0.318740);
		rgb = (0.577060, 0.204783, 0.315654);
		rgb = (0.577596, 0.208664, 0.312565);
		rgb = (0.578135, 0.212545, 0.309542);
		rgb = (0.578676, 0.216431, 0.306516);
		rgb = (0.579214, 0.220287, 0.303496);
		rgb = (0.579746, 0.224106, 0.300518);
		rgb = (0.580271, 0.227977, 0.297566);
		rgb = (0.580793, 0.231817, 0.294618);
		rgb = (0.581315, 0.235646, 0.291715);
		rgb = (0.581835, 0.239463, 0.288810);
		rgb = (0.582353, 0.243268, 0.285910);
		rgb = (0.582870, 0.247097, 0.283066);
		rgb = (0.583386, 0.250916, 0.280201);
		rgb = (0.583901, 0.254739, 0.277381);
		rgb = (0.584416, 0.258531, 0.274552);
		rgb = (0.584931, 0.262342, 0.271740);
		rgb = (0.585443, 0.266156, 0.268980);
		rgb = (0.585951, 0.269966, 0.266198);
		rgb = (0.586456, 0.273771, 0.263439);
		rgb = (0.586961, 0.277575, 0.260676);
		rgb = (0.587466, 0.281374, 0.257925);
		rgb = (0.587972, 0.285180, 0.255221);
		rgb = (0.588478, 0.289013, 0.252494);
		rgb = (0.588984, 0.292818, 0.249767);
		rgb = (0.589491, 0.296652, 0.247081);
		rgb = (0.589999, 0.300465, 0.244376);
		rgb = (0.590507, 0.304300, 0.241716);
		rgb = (0.591016, 0.308135, 0.239031);
		rgb = (0.591526, 0.311969, 0.236379);
		rgb = (0.592038, 0.315846, 0.233692);
		rgb = (0.592548, 0.319698, 0.231058);
		rgb = (0.593055, 0.323559, 0.228420);
		rgb = (0.593562, 0.327429, 0.225773);
		rgb = (0.594071, 0.331309, 0.223134);
		rgb = (0.594583, 0.335229, 0.220510);
		rgb = (0.595095, 0.339131, 0.217865);
		rgb = (0.595609, 0.343048, 0.215226);
		rgb = (0.596126, 0.346976, 0.212613);
		rgb = (0.596645, 0.350921, 0.209994);
		rgb = (0.597164, 0.354880, 0.207388);
		rgb = (0.597680, 0.358830, 0.204776);
		rgb = (0.598196, 0.362821, 0.202147);
		rgb = (0.598721, 0.366829, 0.199533);
		rgb = (0.599248, 0.370837, 0.196964);
		rgb = (0.599771, 0.374879, 0.194370);
		rgb = (0.600294, 0.378931, 0.191738);
		rgb = (0.600819, 0.383009, 0.189149);
		rgb = (0.601346, 0.387090, 0.186548);
		rgb = (0.601874, 0.391215, 0.183949);
		rgb = (0.602403, 0.395345, 0.181345);
		rgb = (0.602933, 0.399486, 0.178782);
		rgb = (0.603464, 0.403678, 0.176158);
		rgb = (0.603995, 0.407873, 0.173594);
		rgb = (0.604521, 0.412102, 0.171015);
		rgb = (0.605043, 0.416348, 0.168436);
		rgb = (0.605562, 0.420618, 0.165848);
		rgb = (0.606084, 0.424928, 0.163317);
		rgb = (0.606609, 0.429252, 0.160731);
		rgb = (0.607129, 0.433600, 0.158195);
		rgb = (0.607639, 0.437998, 0.155649);
		rgb = (0.608144, 0.442412, 0.153086);
		rgb = (0.608644, 0.446848, 0.150582);
		rgb = (0.609134, 0.451324, 0.148071);
		rgb = (0.609610, 0.455826, 0.145615);
		rgb = (0.610079, 0.460356, 0.143119);
		rgb = (0.610542, 0.464933, 0.140685);
		rgb = (0.610991, 0.469544, 0.138267);
		rgb = (0.611421, 0.474170, 0.135829);
		rgb = (0.611833, 0.478839, 0.133514);
		rgb = (0.612226, 0.483539, 0.131212);
		rgb = (0.612600, 0.488287, 0.128920);
		rgb = (0.612950, 0.493049, 0.126718);
		rgb = (0.613275, 0.497875, 0.124574);
		rgb = (0.613572, 0.502705, 0.122487);
		rgb = (0.613837, 0.507592, 0.120512);
		rgb = (0.614069, 0.512502, 0.118669);
		rgb = (0.614264, 0.517459, 0.116848);
		rgb = (0.614418, 0.522434, 0.115160);
		rgb = (0.614530, 0.527456, 0.113657);
		rgb = (0.614594, 0.532510, 0.112266);
		rgb = (0.614607, 0.537595, 0.111032);
		rgb = (0.614566, 0.542708, 0.109999);
		rgb = (0.614468, 0.547849, 0.109114);
		rgb = (0.614308, 0.553016, 0.108421);
		rgb = (0.614082, 0.558212, 0.108010);
		rgb = (0.613787, 0.563446, 0.107850);
		rgb = (0.613419, 0.568682, 0.107943);
		rgb = (0.612974, 0.573946, 0.108312);
		rgb = (0.612449, 0.579232, 0.109026);
		rgb = (0.611842, 0.584522, 0.110040);
		rgb = (0.611148, 0.589820, 0.111320);
		rgb = (0.610353, 0.595132, 0.112963);
		rgb = (0.609471, 0.600443, 0.114856);
		rgb = (0.608494, 0.605748, 0.117169);
		rgb = (0.607411, 0.611060, 0.119811);
		rgb = (0.606215, 0.616350, 0.122763);
		rgb = (0.604930, 0.621618, 0.126124);
		rgb = (0.603536, 0.626876, 0.129757);
		rgb = (0.602026, 0.632107, 0.133692);
		rgb = (0.600413, 0.637306, 0.137967);
		rgb = (0.598689, 0.642469, 0.142496);
		rgb = (0.596862, 0.647588, 0.147334);
		rgb = (0.594916, 0.652662, 0.152416);
		rgb = (0.592872, 0.657697, 0.157790);
		rgb = (0.590707, 0.662667, 0.163419);
		rgb = (0.588441, 0.667579, 0.169258);
		rgb = (0.586085, 0.672429, 0.175280);
		rgb = (0.583613, 0.677213, 0.181507);
		rgb = (0.581049, 0.681916, 0.187985);
		rgb = (0.578388, 0.686560, 0.194586);
		rgb = (0.575646, 0.691121, 0.201310);
		rgb = (0.572809, 0.695614, 0.208243);
		rgb = (0.569878, 0.700018, 0.215285);
		rgb = (0.566888, 0.704346, 0.222470);
		rgb = (0.563814, 0.708597, 0.229738);
		rgb = (0.560662, 0.712753, 0.237171);
		rgb = (0.557458, 0.716845, 0.244622);
		rgb = (0.554182, 0.720839, 0.252219);
		rgb = (0.550853, 0.724766, 0.259874);
		rgb = (0.547470, 0.728605, 0.267574);
		rgb = (0.544043, 0.732376, 0.275394);
		rgb = (0.540571, 0.736058, 0.283238);
		rgb = (0.537067, 0.739685, 0.291141);
		rgb = (0.533507, 0.743228, 0.299094);
		rgb = (0.529936, 0.746702, 0.307079);
		rgb = (0.526333, 0.750112, 0.315113);
		rgb = (0.522696, 0.753461, 0.323192);
		rgb = (0.519049, 0.756752, 0.331281);
		rgb = (0.515367, 0.759983, 0.339437);
		rgb = (0.511681, 0.763162, 0.347595);
		rgb = (0.507990, 0.766293, 0.355785);
		rgb = (0.504280, 0.769372, 0.363984);
		rgb = (0.500550, 0.772410, 0.372217);
		rgb = (0.496820, 0.775405, 0.380485);
		rgb = (0.493085, 0.778365, 0.388763);
		rgb = (0.489350, 0.781287, 0.397049);
		rgb = (0.485614, 0.784180, 0.405376);
		rgb = (0.481884, 0.787038, 0.413711);
		rgb = (0.478142, 0.789866, 0.422057);
		rgb = (0.474411, 0.792674, 0.430440);
		rgb = (0.470680, 0.795455, 0.438824);
		rgb = (0.466955, 0.798219, 0.447235);
		rgb = (0.463220, 0.800964, 0.455667);
		rgb = (0.459518, 0.803693, 0.464121);
		rgb = (0.455810, 0.806409, 0.472577);
		rgb = (0.452124, 0.809110, 0.481054);
		rgb = (0.448436, 0.811796, 0.489555);
		rgb = (0.444772, 0.814472, 0.498091);
		rgb = (0.441108, 0.817144, 0.506616);
		rgb = (0.437487, 0.819803, 0.515175);
		rgb = (0.433858, 0.822465, 0.523755);
		rgb = (0.430280, 0.825110, 0.532352);
		rgb = (0.426720, 0.827756, 0.540960);
		rgb = (0.423186, 0.830401, 0.549598);
		rgb = (0.419708, 0.833036, 0.558241);
		rgb = (0.416257, 0.835673, 0.566923);
		rgb = (0.412868, 0.838305, 0.575612);
		rgb = (0.409520, 0.840937, 0.584314);
		rgb = (0.406245, 0.843562, 0.593044);
		rgb = (0.403035, 0.846190, 0.601780);
		rgb = (0.399905, 0.848819, 0.610541);
		rgb = (0.396872, 0.851439, 0.619320);
		rgb = (0.393950, 0.854061, 0.628104);
		rgb = (0.391152, 0.856683, 0.636905);
		rgb = (0.388472, 0.859301, 0.645709);
		rgb = (0.385935, 0.861918, 0.654530);
		rgb = (0.383585, 0.864526, 0.663367);
		rgb = (0.381407, 0.867128, 0.672196);
		rgb = (0.379424, 0.869728, 0.681023);
		rgb = (0.377672, 0.872325, 0.689863);
		rgb = (0.376170, 0.874907, 0.698686);
		rgb = (0.374923, 0.877482, 0.707507);
		rgb = (0.373981, 0.880045, 0.716318);
		rgb = (0.373340, 0.882596, 0.725106);
		rgb = (0.373043, 0.885136, 0.733865);
		rgb = (0.373112, 0.887654, 0.742601);
		rgb = (0.373570, 0.890156, 0.751300);
		rgb = (0.374439, 0.892639, 0.759946);
		rgb = (0.375723, 0.895095, 0.768546);
		rgb = (0.377467, 0.897524, 0.777098);
		rgb = (0.379671, 0.899923, 0.785572);
		rgb = (0.382352, 0.902288, 0.793974);
		rgb = (0.385527, 0.904619, 0.802283);
		rgb = (0.389213, 0.906913, 0.810503);
		rgb = (0.393385, 0.909161, 0.818619);
		rgb = (0.398074, 0.911369, 0.826627);
		rgb = (0.403255, 0.913528, 0.834507);
		rgb = (0.408926, 0.915628, 0.842255);
		rgb = (0.415083, 0.917688, 0.849859);
		rgb = (0.421704, 0.919678, 0.857309);
		rgb = (0.428791, 0.921615, 0.864606);
		rgb = (0.436305, 0.923489, 0.871734);
		rgb = (0.444231, 0.925293, 0.878682);
		rgb = (0.452541, 0.927032, 0.885454);
		rgb = (0.461203, 0.928705, 0.892037);
		rgb = (0.470211, 0.930311, 0.898424);
		rgb = (0.479521, 0.931839, 0.904620);
		rgb = (0.489103, 0.933297, 0.910617);
		rgb = (0.498950, 0.934685, 0.916408);
		rgb = (0.509019, 0.936004, 0.922005);
		rgb = (0.519281, 0.937246, 0.927394);
		rgb = (0.529715, 0.938416, 0.932588);
		rgb = (0.540292, 0.939517, 0.937592);
		rgb = (0.550997, 0.940549, 0.942401);
		rgb = (0.561804, 0.941509, 0.947020);
		rgb = (0.572686, 0.942411, 0.951459);
		rgb = (0.583621, 0.943243, 0.955728);
		rgb = (0.594606, 0.944015, 0.959825);
		rgb = (0.605610, 0.944731, 0.963765);
		rgb = (0.616637, 0.945388, 0.967563);
		rgb = (0.627648, 0.945989, 0.971214);
		rgb = (0.638645, 0.946543, 0.974739);
		rgb = (0.649620, 0.947052, 0.978146);
		rgb = (0.660548, 0.947515, 0.981449);
		rgb = (0.671439, 0.947934, 0.984653);
		rgb = (0.682276, 0.948316, 0.987765);
		rgb = (0.693064, 0.948662, 0.990803);
		rgb = (0.703779, 0.948977, 0.993775);
}}

\NewDocumentCommand{\TT}{o}{
    \IfNoValueTF {#1} {%
        \vecop T
    }
    {
        \vecop T_{#1}
    }
}

\NewDocumentCommand{\TA}{o}{
    \IfNoValueTF {#1} {%
        \vecop T_{\kern-2pt\mr{A}}
    }
    {
        \vecop T_{\kern-2pt\mr{A},#1}
    }
}

\NewDocumentCommand{\V}{o}{
    \IfNoValueTF {#1} {%
        \op V
    }
    {
        \op V_{#1}
    }
}

\NewDocumentCommand{\W}{o}{
    \IfNoValueTF {#1} {%
        \op W
    }
    {
        \op W_{#1}
    }
}

\NewDocumentCommand{\TPhi}{o}{
    \IfNoValueTF {#1} {%
        \vecop T_{\kern-2pt\Phiup}
    }
    {
        \vecop T_{\kern-2pt\Phiup,#1}
    }
}

\NewDocumentCommand{\matTA}{o}{
    \IfNoValueTF {#1} {%
        \mat T_\mr{A}   
        }
    {
        \mat T_{\mr{A},#1}
    }
}

\NewDocumentCommand{\matTPhi}{o}{
    \IfNoValueTF {#1} {%
        \mat T_\Phiup
        }
    {
        \mat T_{\Phiup,#1}
    }
}

\NewDocumentCommand{\matV}{o}{
    \IfNoValueTF {#1} {%
        \mat V   
        }
    {
        \mat V_{#1}
    }
}

\newcommand{\fk}{f\kern-2pt} 
\usepackage[style=ieee, backend=biber, isbn=false, maxbibnames=9]{biblatex}

\usepackage{nicefrac}
\usepackage{yfonts}
\usepackage{comment}
\graphicspath{{figures/}}
\DeclareUnicodeCharacter{2212}{\ensuremath{-}}

\newcommand{\edits}[1]{}

\begin{document}

	\title{A Stabilized Multilevel B-Spline-Based Fast Integral Method for the Solution of the Electric Field Integral Equation}

	%

	\author{Danijel Jukić,~\IEEEmembership{Graduate Student Member,~IEEE,} Bernd Hofmann,~\IEEEmembership{Member,~IEEE,} Thomas F. Eibert,~\IEEEmembership{Senior Member,~IEEE,} and~Simon  B. Adrian,~\IEEEmembership{Senior Member,~IEEE}
		\thanks{Manuscript received Month xx, 2026; revised Month xx, 20XX.}%
        \thanks{This work was funded by the Deutsche Forschungsgemeinschaft (DFG, German Research Foundation) ‒ 504345461.}%
		\thanks{D. Jukić and S. B. Adrian are with the Fakultät für Informatik und Elektrotechnik, Universität Rostock, 18059 Rostock, Germany (danijel.jukic@uni-rostock.de).}%
        \thanks{Bernd Hofmann is with the Department of Electrical Engineering, Stanford University, Stanford, CA 94305, USA}%
        \thanks{T. F. Eibert is with the Department of Electrical Engineering, School of Computation, Information and Technology, Technical University of Munich, 80290 Munich, Germany}%
        \thanks{Digital Object Identifier}%
	}

	%
	%

\markboth{}%
{Jukić \MakeLowercase{\textit{et al.}}: A Stabilized Multilevel B-Spline-Based Fast Integral Method for the Solution of the Electric Field Integral Equation}
%



\maketitle

\begin{abstract}
    We present a multilevel B-spline-based fast integral method for the solution of the \ac{EFIE}, combining \ac{FFT}-compatible kernel interpolation with robust high-order interpolation.
    Existing \ac{FFT}-accelerated global Lagrange-based approaches rely on equidistant interpolation points and can, therefore, suffer from Runge-type instabilities at high interpolation orders, limiting robust high-accuracy compression.
    In contrast, B-splines on equidistant knot vectors overcome these instabilities and enable robust high-order interpolation for accurate matrix compression.
    Replacing Lagrange interpolation by B-spline interpolation is, however, non-trivial: B-spline coefficients do not coincide with function values at the interpolation points, and the associated sampling matrices can become ill-conditioned.
    To address these challenges, we introduce a knot-removal stabilization strategy, combined with exact interlevel transfers based on knot insertion, yielding accurate, well-conditioned multilevel interpolation.
    Moreover, we propose a factorization strategy that preserves the null space of the scalar potential operator up to machine precision and is compatible with low-frequency preconditioning techniques.
    Numerical results for both canonical and realistic geometries demonstrate robust high-order interpolation without the breakdown observed for Lagrange-based approaches and confirm $\mc O(N)$ complexity.
\end{abstract}
\acresetall


\acresetall
%
\IEEEpeerreviewmaketitle

\section{Introduction}\label{sec:introduction}

\IEEEPARstart{I}{ntegral} equations such as the \ac{EFIE} are a popular choice for the solution of electromagnetic scattering and radiation problems.
They require only a surface discretization of homogeneous scatterers, while the radiation condition is satisfied inherently.
However, a standard \ac{BEM} discretization leads to $\mathcal{O}(N^2)$ complexity for the \ac{MVP}, the assembly time, and the storage requirements with respect to the number of unknowns $N$, which severely limits the problem size.
To address this limitation, so-called fast integral methods have been developed.
These methods typically reduce computational complexity by representing the kernel with a separable expansion, enabling efficient computation of interactions between well-separated basis functions.
For high-frequency scattering problems, fast integral methods such as the \ac{MLFMM}~\cite{coifmanFastMultipoleMethod1993,songMultilevelFastMultipole1997,chewFastEfficientAlgorithms2001} are widely used.
Alternative approaches include butterfly methods~\cite{michielssenMultilevelMatrixDecomposition1996,guoButterflyBasedDirectIntegralEquation2017}, \ac{FFT}-accelerated methods such as the \ac{AIM}~\cite{bleszynskiAIMAdaptiveIntegral1996,yangThreeDimensionalAdaptiveIntegral2012,weiMoreScalableEfficient2014}, and hybrid \ac{MLFMM}-\ac{FFT} schemes~\cite{taboadaHighScalabilityFMMFFT2009,taboadaMLFMAFFTParallelAlgorithm2013}.

In the low-frequency regime, algebraic techniques such as the \ac{ACA}~\cite{zhaoAdaptiveCrossApproximation2005} and its nested variant~\cite{bebendorfConstructingNestedBases2012,zhaoFastNestedCross2019,tetznerIncompleteAdaptiveCross2026} provide efficient compression, but generally lack rigorous error control.
Further approaches include panel clustering~\cite{hackbuschFastMatrixMultiplication1989}, \ac{FFT}-based methods such as the precorrected \ac{FFT} method~\cite{phillipsPrecorrectedFFTMethodElectrostatic1997}, in which the precorrection stage incurs an additional computational cost, and the \ac{ACE}~\cite{shankerAcceleratedCartesianExpansions2007}.
Moreover, variants of the \ac{MLFMM} have been extended to the low-frequency regime by pushing the $k$-space representations into the complex plane~\cite{jiangLowFrequencyFast2004,wulfEfficientImplementationCombined2009a,bogaertNondirectivePlaneWave2008,bogaertLowFrequencyStable2009}.

A different class of fast integral methods is based on kernel interpolation, whose main advantages are kernel independence and controllable error: only kernel samples are required, and the framework introduced in~\cite{hackbuschH2matrixApproximationIntegral2002,bormApproximationIntegralOperators2005} provides theoretical error bounds.
The corresponding kernel expansions are constructed using Lagrange interpolation.
Their application to the \ac{EFIE} was first presented in~\cite{wangNewMultilevelGreens2005}, where translation operators are additionally compressed via a QR-decomposition, introducing additional errors and directional dependence.
Extensions using radial basis functions have been proposed to address electrically larger problems, that is, scatterers whose size exceeds one wavelength~\cite{wangImplementationMultilevelGreens2007,yanshiComparisonInterpolatingFunctions2010,zhaoImprovedFullWaveMultilevel2017}, but they do not allow for exact interlevel translations.

A multilevel Lagrange-based interpolation method has been combined with \ac{FFT} acceleration and applied to the \ac{EFIE} in~\cite{schobertLowFrequencySurfaceIntegral2012,schobertFastIntegralEquation2012} by exploiting the three-level Toeplitz structure of the intralevel translation matrices.
This reduces the \ac{MVP} complexity from $\mathcal{O}(p^6 N)$ to $\mathcal{O}(p^3 \log (p) N)$, where $p$ denotes the polynomial degree of the Lagrange interpolation, enabling higher accuracy and larger problem sizes.
The resulting method is related to precorrected \ac{FFT} approaches~\cite{phillipsPrecorrectedFFTMethodElectrostatic1997}, but employs a multilevel formulation without a precorrection step.
However, \ac{FFT} acceleration requires equidistant interpolation points, which can lead to numerical instabilities due to Runge's phenomenon.
As a consequence, the interpolation error can no longer be reliably controlled by increasing the polynomial degree.

To overcome the instability of high-order equidistant Lagrange interpolation, we propose a B-spline-based kernel interpolation framework for the \ac{EFIE}.
In contrast to global Lagrange interpolation, B-splines mitigate Runge-type instabilities while remaining compatible with equidistant sampling, thereby enabling \ac{FFT} acceleration without the breakdown observed for high-order Lagrange interpolation.
However, B-spline interpolation coefficients do not coincide with function values at the interpolation points.
Consequently, B-spline interpolation requires associated sampling matrices that can become increasingly ill-conditioned at higher polynomial degrees, limiting the stability of high-order interpolation in the multilevel fast integral method.
We overcome this issue by introducing a stabilization strategy that selectively removes knots while preserving exact interlevel transfers via knot insertion.
This yields an accurate and well-conditioned construction of translation and transfer operators.
Critically, we maintain the $\mathcal{O}(\numSplines^3 \log (\numSplines) N)$ \ac{MVP} complexity by leveraging the compact support of B-splines, which leads to sparse interlevel transfers and banded sampling matrices, whose inversion scales linearly with the number of splines $\numSplines$.
In addition, we introduce a star-based factorization strategy for the \ac{EFIE} that preserves the null-space of the scalar potential operator to machine precision and is therefore compatible with low-frequency stabilization and preconditioning techniques.
Numerical results demonstrate stable high-order interpolation without the breakdown observed for Lagrange-based approaches and confirm $\mc O(N)$ scaling with respect to the number of unknowns $N$ for both canonical and realistic geometries.
Preliminary results have been presented at conferences~\cite{jukicInvestigatingBSplineBased2025, jukicBSplineH2,jukicBroadbandStabilizedMultilevel2025a}.

This article is organized as follows:
\Cref{sec:Notation} introduces the \ac{EFIE} and establishes the notation used throughout the paper.
\Cref{sec:singlelevel} presents the single-level algorithm, while \Cref{sec:multilevel} extends it to a multilevel formulation.
\Cref{sec:errorestimates} derives error estimates for the proposed B-spline interpolation and introduces the knot-removal stabilization strategy.
\Cref{sec:factorization} discusses three factorization strategies of the discretized \ac{EFIE}.
Finally, numerical results are presented in \Cref{sec:numerics}.

\section{Background and Notation}\label{sec:Notation}

We consider a time-harmonic incident electromagnetic wave $(\veg e^\mr{i}, \veg h^\mr{i})$ impinging on a \ac{PEC} scatterer with surface $\Gamma$.
The induced surface current density $\veg j$ (normalized by the wave impedance) is related to the incident electric field by the \ac{EFIE}~\cite{mautzHFieldEFieldCombined1977}
\begin{equation}
\left(\TT \veg j \right)_\mr{tan} = (\veg e^\mr{i})_\mr{tan} \,,
\label{eq:efie}
\end{equation}
where the subscript $\mr{tan}$ denotes the tangential component on $\Gamma$.
The \ac{EFIE} operator $\TT$ at wavenumber $k$ is given by
\begin{equation}
\TT \veg j = \jm k\, \TA \veg j - (\jm k)^{-1} \TPhi \veg j \,,
\end{equation}
with the vector potential operator
\begin{align}
    \nonumber \\[-8mm]
    (\TA \veg j)(\veg x) &= \int_\Gamma g(\veg x, \veg y) \veg j (\veg y)\, \dd S(\veg y)\,, \\[-1mm]
\intertext{and the scalar potential operator} \nonumber\\[-7mm]
    (\TPhi \veg j)(\veg x) &= \grad \int_\Gamma g(\veg x, \veg y) \Div_\Gamma \! \veg j (\veg y) \,\dd S(\veg y)\,,
\end{align}
where
\begin{equation}
g(\veg x, \veg y) = \frac{\e^{-\jm k \abs{\veg x - \veg y}}}{4 \uppi \abs{\veg x - \veg y}}
\end{equation}
is the free-space Green's function and $\jm$ denotes the imaginary unit.
Note that a time dependency of $\e^{\,\jm \omega t}$ is assumed and suppressed.

Following a Galerkin approach employing \ac{RWG} functions~\cite{raoElectromagneticScatteringSurfaces1982}, the unknown current density is approximated as $\veg j \approx \sum_{n=1}^N \sbr{\vec j}_n \veg f_{\!n}$, where in contrast to the standard definition, we omit the edge-length normalization in the \ac{RWG} definition as in~\cite[Eq.~14]{adrianElectromagneticIntegralEquations2021}.
This leads to the discretized \ac{EFIE}
\begin{equation}
    \mat T \vec j  = ( \jm k \matTA  - (\jm k)^{-1}\matTPhi ) \vec j = \vec e^\mr{i}\, , \label{eq:discretizedEFIE}
\end{equation}
where
\begin{align}
    \sbr{\matTA}_{mn}\! &= \int_{\Gamma} \!\! \int_{\Gamma} \! \veg f_{\!m}(\veg x) \cdot g(\veg x, \veg y) \veg f_{\!n}(\veg y) \, \dd S(\veg y) \dd S(\veg x) \label{eq:matTA}\,, \\[-0mm]
\intertext{and} \nonumber\\[-6mm]
    \sbr{\matTPhi}_{mn} \! &= -\int_{\Gamma} \! \! \int_{\Gamma} \! \Div_\Gamma{\!\veg f_{\!m}(\veg x)} g(\veg x, \veg y) \Div_\Gamma {\!\veg f_{\!n}(\veg y)} \,\dd S(\veg y) \dd S(\veg x) \label{eq:matTPhi}\,,  \nonumber\\[-2mm]
\end{align}
\\[-6mm]
and the right-hand side
\begin{align}
    [\vec{e}^\mr{i}]_m = \int_{\Gamma} \veg f_{\!m} (\veg y) \cdot \veg e^\mr{i} (\veg y) \,\dd S(\veg y)\,.
\end{align}
After solving \eqref{eq:discretizedEFIE}, for instance, using an iterative solver, the scattered fields are obtained by evaluating the corresponding radiation potentials.
The matrices $\matTA$ and $\matTPhi$ are dense, resulting in $\mc O(N^2)$ cost for assembly, storage, and \ac{MVP}, which becomes prohibitive for many realistic problems.
In the next sections, we establish a robust algorithm that reduces the cost to $\mc O(N)$ (i.e., when the frequency and the scatterer's size are kept fixed).

\section{The Single-Level Acceleration Algorithm} \label{sec:singlelevel}
In this section, we introduce the single-level algorithm, which serves as the building block for the multilevel formulation presented in~\Cref{sec:multilevel}.
Fast integral methods typically reduce the computational complexity by approximating the Green's function through separable expansions, enabling efficient evaluation of interactions between well-separated basis functions.
In the proposed approach, this separation is achieved via polynomial interpolation.
Specifically, we employ a box-based interpolation domain, where test and trial functions are associated with axis-aligned boxes.
For the single-level algorithm, the number of boxes in each spatial direction may differ.
In the multilevel case, these boxes are organized hierarchically using an octree.

\subsection{Polynomial Factorization of the Green's Function} \label{ssec:factorization}
To factorize the Green's function, we employ three-dimensional tensor-product polynomials constructed from one-dimensional polynomials $\ell_i$ along each Cartesian axis $x_i, i\in \{1,2,3\}$~\cite{hackbuschH2matrixApproximationIntegral2002}.
Specifically, we define
\begin{equation}
    \ell_{i} (\veg x) = \ell_{i_1}\!(x_1) \ell_{i_2}\!(x_2) \ell_{i_3}\!(x_3) ,\quad i = 1,\dots,\numSplines^3 \,,
\end{equation}
where $\numSplines$ is the number of polynomials in one Cartesian direction, and the mapping between indices $i$ and $(i_1,i_2,i_3)$ can be chosen arbitrarily.
These polynomials are scaled and shifted to match the domain of the axis-aligned boxes.
For a box centered at $\veg c_{\!\veg x}$, the corresponding scaled and shifted polynomial is denoted by $\ell_{i, \veg c_{\!\veg x}}\!(\veg x)$.

\begin{figure}[tp]
    \centering
    \makebox[\columnwidth][c]{%
        \subfloat[]{
            \includegraphics{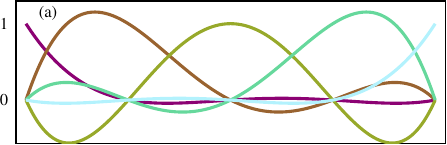}
            \label{fig:lagrange_single_level}
        }
    }\\[-0.4cm]
    \makebox[\columnwidth][c]{%
        \subfloat[]{
            \includegraphics{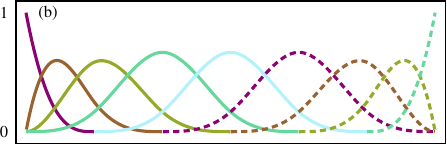}
            \label{fig:bsplines_3_9_single_level}
        }
    }\\[-0.4cm]
    \caption{One-dimensional polynomials $\ell_i(x)$ used for kernel interpolation: (a) Lagrange polynomials of order $p=4$ defined at equispaced interpolation points and (b) cubic B-splines with $\numSplines=9$ defined on an open, equispaced knot vector.}
    \label{fig:lagrange_bspline_single_level}
\end{figure}

For the one-dimensional polynomials $\ell_i(x)$, we consider two choices.
The first choice employs Lagrange polynomials $l_i (x)$ of order $p=s-1$~\cite{berrutBarycentricLagrangeInterpolation2004}, defined at equispaced interpolation points $x_i$ (see \Cref{fig:lagrange_single_level}), representing the state-of-the-art interpolation approach.
The second choice employs $\numSplines$ B-splines $b_i^p\!(x)$ of degree $p$, defined on an open and equispaced knot vector~\cite{deboorPracticalGuideSplines2001} with elements
\begin{equation}
    u_i = \begin{cases}
        0\,, &i \in \sbr{1, p+1} \,,\\
        \frac{i-p-1}{\numSplines-p}\,, &i \in \sbr{p+2, \numSplines}\,,\\
        1\,, &i \in \sbr{\numSplines+1, \numSplines+p+1} \,.\\
    \end{cases} \label{eq:eqknots}
\end{equation}
These B-splines form the foundation for the proposed interpolation framework.
An example illustrating the compact support and overlap of the B-splines is shown in \Cref{fig:bsplines_3_9_single_level}.

\begin{figure}
    \centering
    \includegraphics[width=0.75\linewidth]{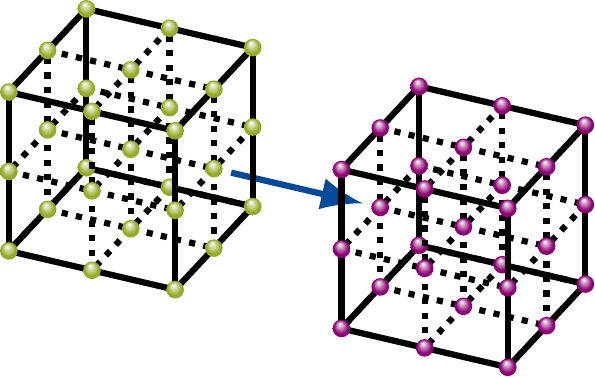}
    \caption{Interpolation grid associated with the axis-aligned boxes of the test and trial functions, shown for three interpolation points along each Cartesian axis.}
    \label{fig:knots_translations}
\end{figure}

Lagrange polynomials are global basis functions with support over the entire interpolation domain and are controlled solely by the polynomial degree $p$.
In contrast, B-splines have local support and are parameterized by both the polynomial degree $p$ and the number of splines $\numSplines$, allowing the spline degree and the number of splines to be controlled independently.
The interpolation points $x_i$ are chosen equispaced (to later allow for \ac{FFT} acceleration).
\Cref{fig:knots_translations} exemplifies the placement of the interpolation points relative to the axis-aligned boxes of a set of test and trial functions.

To illustrate how polynomial-based interpolation is used to factorize the Green's function $g(\veg x, \veg y)$, consider the discretized vector potential
\begin{multline}
    \sbr{\matTA}_{mn} \!\approx \sum_{i=1}^{\numSplines^3} \int_{\supp \veg f_{\!m}} \veg f_{\!m}(\veg x) \ell_{i,\veg c_{\!\veg x}}\! (\veg x) \, \dd \veg x \cdot \sum_{j=1}^{\numSplines^3} w_{i,j}
    \\ \int_{\supp \veg f_{\!n}} \veg f_{\!n}(\veg y) \ell_{j,\veg c_{\!\veg y}}\! (\veg y) \, \dd S(\veg y) \eqcolon \sum_{d=1}^3 \vec m_{\,d,\veg f_{\!m}}^\T \mat W \vec m_{\,d,\veg f_{\!n}} \label{eq:TAfactor}  \,,
\end{multline}
where $\sbr{\mat W}_{ij} = w_{i,j}$ denotes the matrix of interpolation weights, which will be specified in the next paragraph.
The right-hand side of \eqref{eq:TAfactor} is obtained by regrouping the terms according to the three Cartesian components of the vector-valued \acp{RWG} $\veg f_{\!n}$.
This allows the integrals, defined over the supports $\supp \veg f_{\!m}$ and $\supp \veg f_{\!n}$, to be expressed as a sum over $d=1,2,3$, where $\vec m_{d,\veg f_{\!m}}$ and $\vec m_{d,\veg f_{\!n}}$ denote the corresponding moment vectors of the test and trial functions for each spatial dimension.

For Lagrange polynomials, the weights $w_{i,j}$ are the samples of $g(\veg x, \veg y)$ at the interpolation points, that is,
\begin{equation}
    \mat W^{\mr{Lagrange}}= \sbr{\mat G}_{m,n} \coloneq g(\veg x_m, \veg y_n)\,.
\end{equation}
In contrast to Lagrange polynomials, the weights of the B-spline interpolation do not correspond to the samples of the Green's function.
Instead, the weights are obtained via~\cite[Chapter~9.2]{pieglNURBSBook1997}
\begin{equation}
    \mat W^{\mr{B-spline}} = {\left(\mat B_\numSplines^{\text{3D},p}\right)}^{-\T} \mat G {\left(\mat B_\numSplines^{\text{3D},p}\right)}^{-1} \,,
\end{equation}
where the sampling matrix $\mat B_\numSplines^{\text{3D},p}$ contains the B-splines evaluated at the interpolation points
\begin{equation}
    \sbr{\mat B_\numSplines^{\text{3D},p}}_{m,n} = b_{m}^p (\veg x_n)\,.
\end{equation}

At first glance, the required inversion of $\mat B_\numSplines^{\text{3D},p}$ appears computationally expensive, suggesting a cost of $\mathcal{O}(\numSplines^6)$.
However, this estimate is overly pessimistic due to the sampling matrix's structure.
First, the tensor-product structure of the B-splines and the interpolation grid implies
\begin{equation}
    \mat B_\numSplines^{\text{3D},p} = \mat B_\numSplines^{\text{1D},p} \otimes \mat B_\numSplines^{\text{1D},p} \otimes \mat B_\numSplines^{\text{1D},p}\,,
\end{equation}
where $\mat B_{\numSplines}^{\text{1D},p}$ denotes the one-dimensional sampling matrix and $\otimes$ the tensor (Kronecker) product, that is,~\cite[Definition~4.2.1]{hornTopicsMatrixAnalysis2010}
\begin{equation}
    \mat A \otimes \mat B =
  \begin{bmatrix}
      \sbr{\mat A}_{1,1} \mat B & \cdots & \sbr{\mat A}_{1,n} \mat B \\
                 \vdots & \ddots &           \vdots \\
      \sbr{\mat A}_{n,1} \mat{B} & \cdots & \sbr{\mat A}_{n,n} \mat B
\end{bmatrix}\,.
\end{equation}
Second, owing to the compact support of B-splines, $\mat B_{\numSplines}^{\text{1D},p}$ is banded with semi-bandwidth $p-1$ and, therefore, its inverse can be applied in $\mathcal{O}(\numSplines)$ operations without explicit construction.
Exploiting the tensor-product structure, applying the inverse of the three-dimensional sampling matrix, thus, reduces to a computational cost of $\mathcal{O}(\numSplines^3)$~\cite[Chapter~9.2]{pieglNURBSBook1997}.

Although $(\mat B_{\numSplines}^{\text{3D},p})^{-1}$ can be applied efficiently to a vector $\vec m$, the result may still suffer from numerical inaccuracies due to the ill-conditioning of $\mat B_{\numSplines}^{\text{3D},p}$.
In particular, the condition number deteriorates with increasing $p$ and $\numSplines$, which can adversely affect the accurate evaluation of $\mat W^{\mr{B-spline}} \vec m_{\,d,\veg f_{\!n}}$.
We defer a detailed discussion of this issue until after introducing B-spline-B-spline interpolation in the context of interlevel transfers, which forms the basis for the stabilization strategy presented in \Cref{ssec:stab}.

\subsection{Protruding Basis Functions}\label{sec:protruding}
The support of the interpolation polynomials is defined with respect to the axis-aligned boxes.
However, the support of the underlying basis functions may extend beyond these boxes.
This must be taken into account when computing the moment vectors to avoid truncation errors.

\begin{figure}[tp]
	\centering
	\subfloat[][a)]{\includegraphics[scale=0.55]{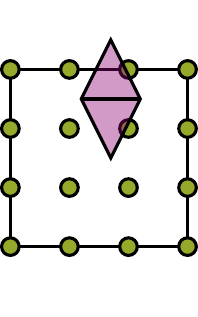}\label{fig:protrusion_original}} 	\hfill
	\subfloat[][b)]{\includegraphics[scale=0.55]{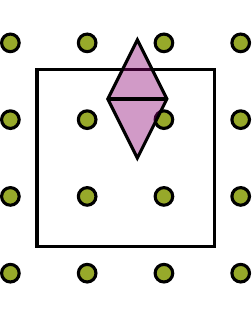}\label{fig:protrusion_scaled}} 	\hfill
	\subfloat[][c)]{\includegraphics[scale=0.55]{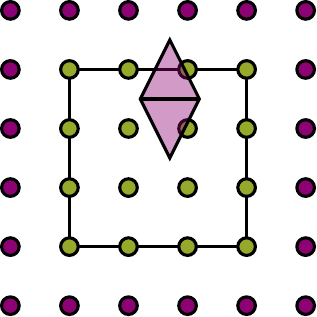}\label{fig:protrusion_buffer}}
	\caption{Interpolation points and support of a protruding \ac{RWG} function. From left to right: (a) original interpolation points, (b) scaled and shifted interpolation points, and (c) interpolation points with additional buffer knots.}
	\label{sphereCurrent}
\end{figure}

For Lagrange polynomials, two strategies have been proposed~\cite{schobertLowFrequencySurfaceIntegral2012,schobertFastIntegralEquation2012}: either extrapolating $\ell_{i, \veg c_{\!\veg x}}\!(\veg x)$ (i.e., evaluating $\ell_{i, \veg c_{\!\veg x}}\!(\veg x)$ outside the box) or rescaling the interpolation points $x_i$ such that the support of the basis function is fully contained.
These approaches are illustrated in \Cref{fig:protrusion_original,fig:protrusion_scaled}.
While rescaling the $x_i$ can mitigate the effects of the Runge phenomenon, it does not eliminate it.
Moreover, in a multilevel setting, extrapolation introduces errors in the interlevel translations.

For B-splines, both strategies are unsuitable.
We have observed that extrapolation performs poorly, as only a few B-splines near the box boundary capture the contributions of protruding basis functions, resulting in limited accuracy.
Rescaling the interpolation points is also not possible, since B-splines are tied to the knot vector.
In particular, such rescaling destroys the uniform structure of the knot vectors across levels, which is required to enable exact interlevel translations (see \Cref{ssec:bsplinetransfer,ssec:constsplines,ssec:moresplines}).

Instead, we introduce a buffer-knot approach, available only in the B-spline approach, that enables exact interlevel translations.
Specifically, we extend the knot vector by adding $\numSplines_\mr{b}$ buffer knots on each side (i.e., increasing the number of splines per axis by $2\numSplines_\mr{b}$), ensuring that the support of every basis function is fully contained within the enlarged interpolation domain, as illustrated in \Cref{fig:protrusion_buffer}.
The required number of buffer splines on each side of the knot vector is
\begin{equation}
    \numSplines_\mr{b} \coloneq \lceil \Updelta_\text{max}^\text{rel} /\Updelta u\rceil\,,
\end{equation}
where $\lceil \cdot \rceil$ denotes the ceiling function, that is, rounding to the smallest integer greater than or equal to its argument,
\begin{equation}
    \Updelta u = \frac{1}{\numSplines-p}
\end{equation}
denotes the distance between adjacent knots for an open, uniformly spaced knot vector, and
$\Updelta_\text{max}^\text{rel}$ denotes the maximum relative protrusion of a basis function support beyond its associated box, taken over all boxes on a given level.
This choice ensures that the enlarged knot vector fully covers the support of all basis functions.

\subsection{FFT Acceleration}
A direct evaluation of \eqref{eq:TAfactor} exhibits a computational complexity of $\mathcal{O}(\numSplines^6)$, which rapidly becomes prohibitive as the number of interpolation polynomials increases.
In the Lagrange case, the interpolation-weight matrix $\mat W^\text{Lagrange} = \mat G$ possesses a three-level Toeplitz structure.
This structure has previously been exploited for efficient \ac{FFT}-accelerated \acp{MVP}~\cite{schobertLowFrequencySurfaceIntegral2012}, reducing the computational cost to $\mathcal{O}(\numSplines^3 \log \numSplines)$.

In contrast, the B-spline interpolation-weight matrix $\mat W^\mr{B-spline}$ does not directly exhibit a Toeplitz structure.
However, by transforming the moment vectors $\vec m_{d, \veg f_{\!m}}$ with the B-spline sampling matrix $\mat B_\numSplines^{\text{3D},p}$, the factorized form in \eqref{eq:TAfactor} can be reformulated as
\begin{equation}
    \sbr{\matTA}_{mn} \!\approx \sum_{d=1}^3 \vec m_{\,d,\veg f_{\!m}}^\T \mat W^\text{B-spline} \vec m_{\,d,\veg f_{\!n}} =  \sum_{d=1}^3 \wvec m_{\,d,\veg f_{\!m}}^\T \mat G \wvec m_{\,d,\veg f_{\!n}} \label{eq:TAfactorBspline}  \,,
\end{equation}
where the transformed moments are defined as
\begin{equation}
    \wvec m_{d, \veg f_{\!m}} = \left(\mat B_\numSplines^{\text{3D},p}\right)^{-1} \vec m_{d, \veg f_{\!m}}\,.
\end{equation}
Since systems involving $\mat B_\numSplines^{\text{3D},p}$ can be solved in $\mathcal{O}(\numSplines^3)$, the overall complexity of the accelerated factorization remains $\mathcal{O}(\numSplines^3 \log \numSplines)$, matching that of the Lagrange-based formulation.

\section{The Multilevel Acceleration Algorithm} \label{sec:multilevel}
While the single-level algorithm already reduces the computational cost of interactions within a fixed level (see \eqref{eq:TAfactor}), its overall complexity remains superlinear in the number of unknowns.
For the Lagrange interpolation, the multilevel algorithm introduces a hierarchical structure that reduces the \ac{MVP}, storage, and assembly complexity to $\mathcal{O}(N)$ for fixed $p$ and $\numSplines$ \cite{hackbuschH2matrixApproximationIntegral2002};
the \ac{FFT} acceleration reduces the complexity from $\mc{O}(\numSplines^6 N)$ to $\mc{O} (\numSplines^3 \log (\numSplines) N)$~\cite{schobertLowFrequencySurfaceIntegral2012}.
This improved complexity is maintained by our algorithm.

To this end, interlevel transfers are incorporated alongside intralevel translations.
The algorithm is formulated on a regular octree hierarchy~\cite{hackbuschH2matrixApproximationIntegral2002} to exploit the translation invariance of the kernel.
As a consequence, the translation matrices $\mat W^\text{Lagrange}$ and $\mat W^\text{B-spline}$ are themselves translation invariant.
Therefore, only a constant number of distinct translation operators must be stored per level, yielding a substantial reduction in memory requirements.

A multilevel formulation additionally requires interlevel translation operators that map moment vectors from child boxes to their parent box.
Due to the linearity of the integrals in \eqref{eq:matTA}, this operation can be realized by interpolating the parent-level polynomials using the child-level polynomials.
This representation is exact when the polynomial spaces are nested.
Specifically, the parent-level moments are obtained as
\begin{equation} \label{eq:interleveltransfer}
    \vec m_{d, \veg f_{\!m}}^\text{parent} = \sum_{c=1}^{\# \,\text{children}} \mat C_c \vec m_{d, \veg f_{\!m}}^{\text{child}_c} \,,
\end{equation}
where we will detail the construction of $\mat C_c$ in the following.

\subsection{Review of Lagrange Interlevel Transfer}
\begin{figure}[tp]
    \centering
    \makebox[\columnwidth][c]{%
        \subfloat[]{
            \includegraphics{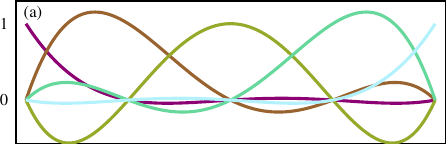}
            \label{fig:lagrange}
        }
    }\\[-0.4cm]
    \makebox[\columnwidth][c]{%
        \subfloat[]{
            \includegraphics{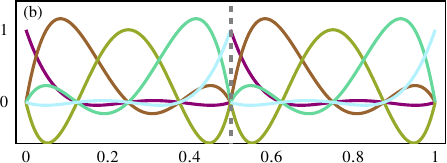}
            \label{fig:lagrange_children}
        }
    }\\[-0.4cm]
    \caption{(a) Set of $\numSplines=5$ parent Lagrange polynomials of order $p=4$ and (b) the two sets of children Lagrange polynomials of order $p=4$.}
    \label{fig:lagrange_interlevel}
\end{figure}

To motivate the choices we make for our new B-spline approach, we first review the Lagrange interpolation case.
In that case, the parent polynomial is expressed in terms of the child polynomials by sampling the parent polynomial at the interpolation points on the child level.
If the polynomial degree is the same on both levels, this representation is exact.
Otherwise, if the degree is increased on the parent level, the parent polynomial cannot be represented exactly, and the transfer corresponds to an interpolation of the parent polynomial.
In one dimension, the corresponding transfer matrices are~\cite{hackbuschH2matrixApproximationIntegral2002}
\begin{equation}
    \sbr{\mat C^{\text{1D}}_c}_{mn} = l_m^\mr{parent} (x^c_n) \,,
\end{equation}
where $x^c_n$ is the $n$th interpolation point of the $c$th child.
A set of parent and child polynomials is depicted in \Cref{fig:lagrange_interlevel}.

In one dimension, only two distinct transfer matrices arise, corresponding to the left and right child transfer matrices $
    \mat C^\text{1D} \in \{\mat C^{\mr{L}}, \mat C^{\mr{R}}\}$.
In three dimensions, due to the regular octree structure, each child box corresponds to one of the $2^3=8$ combinations of  ${C^\mr{L}}$ and ${C^\mr{R}}$.
The corresponding three-dimensional transfer matrices are given by tensor products
\begin{equation}
    \mat C^\text{3D} \in \{ \mat C^{c_1} \otimes \mat C^{c_2} \otimes \mat C^{c_3} \,|\, c_i \in \{\mr{L},\mr{R}\} \}.
\end{equation}
Consequently, only $\mat C^\mr{L}$ and $\mat C^\mr{R}$ need to be stored per level, and all interlevel translations are obtained via the tensor-product structure.

At low frequencies, that is, when the scatterer size is smaller than $0.5\lambda$, $p$ is kept constant across all levels, as is standard in the static regime~\cite{bormApproximationIntegralOperators2005}.
For boxes with edge length greater than or equal to $0.5\lambda$, $p$ is doubled at coarser levels to maintain accuracy~\cite[p.~40]{schobertSchnelleLosungElektromagnetischer}.
In this case, the interlevel transfer is no longer exact, since a higher-degree polynomial cannot be represented exactly by lower-degree polynomials.
A discussion of how this affects the overall accuracy of the fast integral method can be found in~\cite{bormApproximationIntegralOperators2005}.

\subsection{B-Spline Interlevel Transfer} \label{ssec:bsplinetransfer}
For B-splines, the transfer matrices could in principle be obtained by solving a \ac{LSE}, similar to the one used for the kernel interpolation.
However, this approach introduces numerical errors that accumulate across levels, reducing the accuracy of the multilevel scheme.
These errors originate from the ill-conditioning of $\mat B_{\numSplines}^{\text{3D},p}$ (see \Cref{ssec:stab}).

Instead, we propose computing the transfer coefficients via so-called knot insertion, yielding an exact and numerically stable analytical B-spline-B-spline interpolation.
Knot insertion requires the scaled and shifted parent knot vector to be contained in the corresponding child knot vectors; otherwise, knots would have to be removed, which would destroy the exactness of the interpolation.
The required nestedness is ensured for the equidistant configurations discussed in \Cref{ssec:constsplines,ssec:moresplines}.
Moreover, because the child spline configurations are again equidistant, the same knot-insertion procedure can be applied recursively across all levels.

Knot insertion is a standard tool in spline theory (introduced for curve refinement and increasing the degrees of freedom in \ac{CAD})~\cite[Chapter~5.2]{pieglNURBSBook1997}:
Given a knot vector
\begin{equation}
    \{u_1,\dots, u_{s+p+1}\}\,,
\end{equation}
a spline curve with coefficients $w_i$ is defined as
\begin{equation}
    \sum_{i=1}^\numSplines b_i^p\!(x) w_i\,.
\end{equation}
The objective of knot insertion is to construct an equivalent spline representation~\cite[Chapter~5.2]{pieglNURBSBook1997}
\begin{equation}
    \sum_{i=1}^{\numSplines+1} \Bar{b}_i^p\!(x) \Bar{w}_i = \sum_{i=1}^\numSplines b_i^p\!(x) w_i\,,
\end{equation}
defined on an augmented knot vector~\cite[Chapter~5.2]{pieglNURBSBook1997}
\begin{equation}
    \{\Bar{u}_1 = u_1,\dots, \Bar{u}_k = u_k, \Bar{u}_{k+1} = \Bar{u},
    \dots, \bar{u}_{s+p+2} = u_{s+p+1}  \}\, ,
\end{equation}
where the new knot $\bar{u}$ is inserted at position $k+1$ such that $u_k \leq \bar{u} < u_{k+1}$, which determines the interval in which the refinement takes place.
The form of the augmented knot vector is given in \Cref{ssec:constsplines,ssec:moresplines}.
The index $k$ can be determined efficiently using a bisection search, after which the coefficients of the refined spline are computed as~\cite[Chapter~5.2]{pieglNURBSBook1997}
\begin{equation}
    \bar{w}_i = \alpha_i w_i + (1-\alpha_i) w_{i-1} \,,
\end{equation}
with~\cite[Chapter~5.2]{pieglNURBSBook1997}
\begin{equation}
    \alpha_i =
    \begin{cases}
        1, & i \leq k-p\,, \\
        \frac{\Bar{u} - u_i}{u_{i+p} - u_i} & k-p+1 \leq i \leq k\,,\\
        0, & i \geq k+1\,.
    \end{cases}
\end{equation}
Only the coefficients in the interval $i \in [k-p+1,k]$ are modified, reflecting the local support of B-splines.

Based on the knot insertion scheme, the interlevel transfer matrix $\mat C^\text{1D}$ is constructed row by row.
For each parent-level B-spline, indexed by $j$, we define an auxiliary curve with coefficients
\begin{equation}
    w^{(j)}_i = \begin{cases}
        1, & i=j\,, \\
        0, & \text{otherwise}\,.
    \end{cases}
\end{equation}
This curve represents the $j$th parent B-spline.
All knots required to match the child-level knot vectors are then inserted sequentially into this parent curve using the knot-insertion procedure described above.
After completion of the knot insertions, the resulting coefficients $w^{(j)}_i$ correspond to the representation of the parent B-spline in terms of child-level B-splines, and the transfer matrix $\sbr{\mat C}_{ji} = w^{(j)}_i$.
Due to the locality of the knot insertion, $\mat C$ is a sparse (banded) matrix.
Note that this form of interlevel transfer is exact (i.e., close to machine precision, as addressed in more detail in Section~\ref{interlevelStab}).

\subsection{Case of a Constant Number of B-Splines} \label{ssec:constsplines}
\begin{figure}[tp]
    \centering
    \makebox[\columnwidth][c]{%
        \subfloat[]{
            \includegraphics{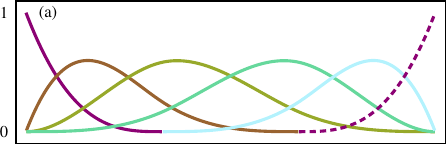}
            \label{fig:bsplines_3_6}
        }
    }\\[-0.4cm]
    \makebox[\columnwidth][c]{%
        \subfloat[]{
            \includegraphics{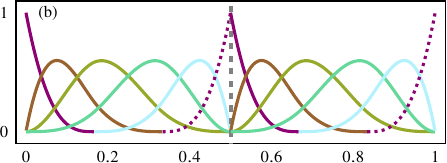}
            \label{fig:bsplines_children}
        }
    }\\[-0.4cm]
    \caption{(a) Set of $\numSplines=6$ parent B-splines of order $p=3$ and (b) the two sets of children B-splines. Notably, the number of B-splines remains constant, with $\numSplines=6$ B-splines in both the parent and child clusters for boxes with an edge length smaller than $0.5\lambda$.}
    \label{fig:bspline_interlevel_lf}
\end{figure}

For parent boxes with an edge length smaller than $0.5\lambda$, the number of splines $\numSplines$ is kept the same on the parent and child levels, as this is sufficient to maintain the desired accuracy.
Consequently, the knot spacing $\Updelta u$ in the parametric domain remains unchanged.
However, after scaling and shifting the knot vectors to the physical boxes, the knot spacing on the parent level is half that on the child level, since the edge length of the child box is half that of the parent box.
A one-dimensional example without buffer knots is shown in \Cref{fig:bspline_interlevel_lf}, illustrating the parent polynomials and the two corresponding sets of child polynomials.

\begin{figure}[tp]
    \centering
    \makebox[\columnwidth][c]{%
        \subfloat[]{
            \includegraphics{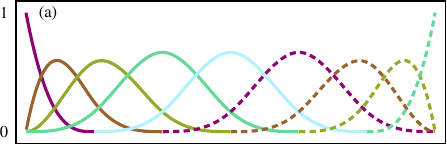}
            \label{fig:bsplines_3_9}
        }
    }\\[-0.4cm]
    \makebox[\columnwidth][c]{%
        \subfloat[]{
            \includegraphics{figures/bsplines_children.pdf}
            \label{fig:bsplines_children_mf}
        }
    }\\[-0.4cm]
    \caption{(a) Set of $\numSplines=9$ parent B-splines of order $p=3$ and (b) the two sets of $\numSplines=6$ children B-splines of order $p=3$. Notably, the number of B-splines increases in the parent cluster for boxes with an edge length greater than or equal to $0.5\lambda$}
    \label{fig:bspline_interlevel_mf}
\end{figure}

Due to the coarser resolution of the scaled knot vector on the parent level, only half as many buffer knots are required as on the child level.
In particular, doubling the physical knot spacing reduces the number of buffer knots required to cover the same relative protrusion by a factor of two.
Thus, for uniform trees, the number of parent-level buffer knots is determined solely by the number of buffer knots on the child levels.
For nonuniform trees, however, leaf boxes on the parent level must also be considered, requiring sufficiently many buffer knots to avoid truncating their contributions.

Combining these requirements, the number of splines and buffer knots on the parent level is given by
\begin{align}
    \numSplines^\text{parent} &= \numSplines^\text{child}\,, \\
    \numSplines_\mr{b}^\text{parent} &= \max \left( \lceil \numSplines_\mr{b}^\text{child}/2 \rceil, \numSplines_\mr{b}^\text{parent, leaves} \right)\,,
\end{align}
where $\numSplines_\mr{b}^\text{parent, leaves}$ denotes the number of buffer knots required to cover the support of parent-level leaf basis functions.

\subsection{Case of an Increasing Number of B-Splines}  \label{ssec:moresplines}
When the edge length of the axis-aligned box exceeds $0.5\lambda$, the number of splines $\numSplines$ is increased on coarser levels to maintain comparable factorization accuracy.
Specifically, the spline spacing in the parametric domain $\Updelta u$ is halved on the parent level compared to the child level.
After scaling and shifting the knot vectors to the physical boxes, the resulting physical knot spacing is identical on both levels.
At the same time, the parent-level knot vector must remain equidistant to ensure that the same construction can be applied recursively on subsequent levels and to enable exact interlevel transfers.
A one-dimensional example without buffer knots is shown in \Cref{fig:bspline_interlevel_mf}.

To obtain a knot vector whose physical knot spacing matches that of the child level, the number of splines on the parent level is
\begin{align}
    \numSplines^\text{parent} &= 2\numSplines^\text{child} - p\,, \\
    \numSplines_\mr{b}^\text{parent} &= \max \left( \numSplines_\mr{b}^\text{child}, \numSplines_\mr{b}^\text{parent, leaves} \right)\,.
\end{align}
As discussed in \Cref{ssec:constsplines}, sufficient buffer knots must be provided on the parent level to fully cover the support of all basis functions.
Since the physical knot spacing on the parent level is identical to that on the child level, the same physical protrusion corresponds to the same number of knot intervals.
Consequently, at least the same number of buffer splines as on the child level is required to avoid truncation.

\subsection{Assessing the Stability of B-Spline Interlevel Transfers} \label{interlevelStab}
\begin{figure}
    \centering
    \includegraphics{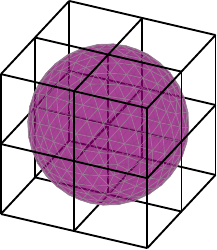}
    \caption{Surface mesh of a sphere with the corresponding axis-aligned boxes of a two-level octree.}
    \label{fig:sphereoctree}
\end{figure}

To assess the numerical stability of the B-spline interlevel transfers, we consider a sphere within a two-level octree.
The surface mesh and corresponding octree boxes are shown in \Cref{fig:sphereoctree}.
The moment vectors on the top level are computed both directly by numerical quadrature and indirectly from the child-level moment vectors using the transfer matrices defined in \eqref{eq:interleveltransfer}.
The number of B-splines $\numSplines$ is kept constant across both levels.

\begin{figure}[tp]
    \centering
    \makebox[\columnwidth][c]{%
        \subfloat[]{
            \includegraphics{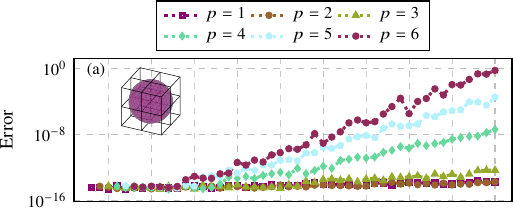}
            \label{fig:bsplines_i2isphere_sub}
        }
    }\\[-0.4cm]
    \makebox[\columnwidth][c]{%
        \subfloat[]{
            \includegraphics{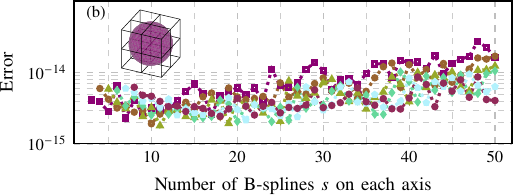}
            \label{fig:bsplines_i2isphere_opt}
        }
    }\\[-0.4cm]
    \caption{Maximum relative error of the B-spline interlevel transfer for increasing numbers of splines. (a) Transfer matrices computed by solving the B-spline sampling \ac{LSE}. (b) Transfer matrices constructed using knot insertion.}
    \label{fig:bsplines_i2i}
\end{figure}

The maximum relative error introduced by the interlevel transfer is shown in \Cref{fig:bsplines_i2i}.
Two constructions of the transfer matrices are compared: a direct construction based on solving a \ac{LSE} using the B-spline sampling matrix $\mat B_\numSplines^{\text{3D},p}$, and the proposed exact approach based on knot insertion as described in \Cref{ssec:bsplinetransfer}.
As shown in \Cref{fig:bsplines_i2i}, the directly constructed transfer matrices exhibit rapidly increasing errors as $\numSplines$ grows, indicating numerical instability.
In contrast, the proposed transfer matrices maintain errors in the order of $\num{1e-14}$, that is, close to machine precision, demonstrating numerical stability even for large $\numSplines$.

\section{Error Estimates and Stability} \label{sec:errorestimates}
We first discuss error estimates for the proposed kernel interpolation scheme and then return to the issue of stable B-spline kernel interpolation for $p>2$, which was deferred in \Cref{ssec:factorization}.
The following analysis illustrates the importance of numerical stability for practical error control.

\subsection{Error Estimates} \label{ssec:error}
In this section, we derive error estimates for tensor-product B-spline interpolation by adapting the framework of~\cite{hackbuschH2matrixApproximationIntegral2002}.
That framework was originally developed for Lagrange interpolation with Chebyshev interpolation points, where the required error and stability estimates are well established.
In contrast, for equidistant interpolation points, which are required to enable the \ac{FFT} acceleration, no such uniformly bounded stability estimates are available.

To derive tensor-product interpolation error estimates, we require an error estimate for the one-dimensional interpolation operator $\Pi$.
For B-splines, this is given by~\cite{schumakerSplineFunctionsBasic2010}
\begin{equation}
    \norm{\Pi f - f}_{\infty,[0,1]} \leq (\Updelta u)^{p+1} \norm{\partial^{p+1}f}_{\infty,[0,1]} \,,
\end{equation}
where $\norm{\cdot}_{\infty,[0,1]}$ denotes the maximum norm on $[0,1]$ and $\Updelta u = 1/(\numSplines-p)$ denotes the maximum knot spacing in the parametric domain.

However, for B-spline interpolation with equidistant knots and equidistant interpolation points, the stability constant associated with the interpolation operator $\Pi$ in estimates of the form
\begin{equation} \label{eq:1dstabilityestimate}
    \norm{\Pi f}_{\infty,[0,1]} \leq  C_2 \norm{f}_{\infty,[0,1]}
\end{equation}
can grow rapidly with the number of splines and the polynomial degree, although it is fixed for a given interpolation setup.
In fact, the growth in $C_2$ is reflected in the increased ill-conditioning of the sampling matrices $\mat B_{\numSplines}^{\text{1D},p}$, which causes numerical errors to accumulate.
As discussed in \Cref{ssec:stab}, we address this issue with a stabilization strategy.

To extend the one-dimensional results to higher dimensions, we consider the tensor-product interpolation operator $\Pi_{B_d}$ defined on a $d$-dimensional box $B_d$ with edge length $h$.
The operator $\Pi_{B_d}$ is obtained as the tensor product of the corresponding one-dimensional interpolation operators in each coordinate direction.

Combining the error and stability estimates yields the $d$-dimensional tensor-product bound~\cite{hackbuschH2matrixApproximationIntegral2002}
\begin{equation} \label{eq:errorapprox}
    \norm{\Pi_{B_d} f - f}_{\infty,B_d} \leq (\Updelta u)^{p+1} C_2^{\,d-1} \left(\frac{h}{2}\right)^{p+1} \sum_{k=1}^{d} \norm{\uppartial_k^{p+1} f}_{\infty,B_d}\,.
\end{equation}
Notably, the tensor-product estimate depends only on the $(p+1)$st derivatives in the coordinate directions and does not involve mixed derivatives.

\subsection{Stabilization of B-Spline Interpolation} \label{ssec:stab}
The relationship between interpolation points $x_i$ and knots $u_i$ plays a critical role in the stability of B-spline interpolation~\cite[pp.~180]{deboorPracticalGuideSplines2001}.
While the Schoenberg-Whitney theorem guarantees uniqueness for equidistant interpolation points and equidistant knots~\cite[pp.~171]{deboorPracticalGuideSplines2001}, the corresponding stability constant $C_2$ in~\eqref{eq:1dstabilityestimate} can grow rapidly with the number of splines and the polynomial degree.
In fact, using~\cite[p.~190]{deboorPracticalGuideSplines2001}, the stability constant $C_2$ in \eqref{eq:1dstabilityestimate} can be bounded by
\begin{equation}
    C_2 \leq \norm{\left(\mat B_{\numSplines}^{\text{1D},p}\right)^{-1}}_\infty
\end{equation}
with the norm
\begin{equation}
    \norm{\mat B}_\infty = \max_m \sum_n \abs{[\mat B]_{m,n}}\,.
\end{equation}
For large numbers of splines $\numSplines$, the resulting sampling matrices $\mat B_{\numSplines}^{\text{1D},p}$ can become severely ill-conditioned, particularly for higher polynomial degrees $p$, substantially reducing interpolation accuracy.

\begin{figure}
    \centering
    \makebox[\columnwidth][r]{%
        \subfloat[]{
            \includegraphics{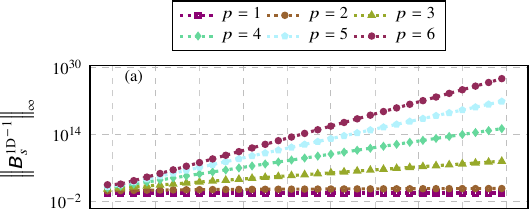}
            \label{fig:cond_equidistant}
        }
    }\\[-0.4cm]
    \makebox[\columnwidth][r]{%
        \subfloat[]{
            \includegraphics{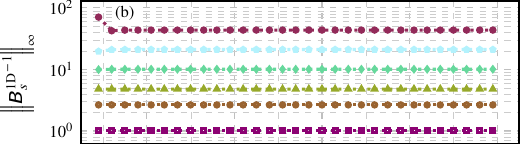}
            \label{fig:cond_avg}
        }
    }\\[-0.4cm]
    \makebox[\columnwidth][r]{%
        \subfloat[]{
            \includegraphics{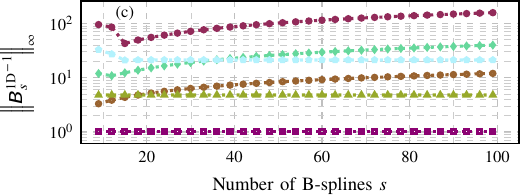}
            \label{fig:cond_stable}
        }
    }\\[-0.4cm]
    \caption{Norm of the inverse B-spline sampling matrix for (a) equidistant knots, (b) averaged knots, and (c) stabilized equidistant knots. Stabilization significantly improves the growth of the stability constant compared to the equidistant case while preserving compatibility with the multilevel scheme.}
    \label{fig:condition_number}
\end{figure}

The growth of the stability constant for equidistant knots and interpolation points is illustrated in \Cref{fig:cond_equidistant}.
A common strategy to improve the interpolation accuracy is the use of averaged knots, obtained by averaging the interpolation points~\cite[Chapter~9.2.1]{pieglNURBSBook1997}.
For equidistant interpolation points, the averaged knot vector is given by
\begin{equation}
    u_i^\text{avg} = \begin{cases}
        0\,, &i \in \sbr{1, p+1} \,,\\
        \frac{2i-p-3}{2(\numSplines-1)}\,, &i \in \sbr{p+2, \numSplines}\,,\\
        1\,, &i \in \sbr{\numSplines+1, \numSplines+p+1} \,.\\
    \end{cases} \label{eq:avgknots}
\end{equation}
These knots significantly reduce the stability constant and thereby improve the interpolation accuracy, as shown in \Cref{fig:cond_avg}.
However, averaged knots do not allow exact interlevel transfers (a requirement for the multilevel algorithm), since not all parent knots are contained within the child knot vectors.

This creates a conflict: averaged knots improve stability, whereas equidistant knots are required for exact interlevel transfers (see \Cref{ssec:bsplinetransfer}).
At the same time, equidistant knots lead to a catastrophic growth in the stability constant.
Therefore, we need a stabilization strategy that preserves the multilevel structure while improving numerical stability.

A key observation of this work is that removing knots near the boundary of an equidistant knot vector significantly improves the stability of the B-spline interpolation.
For $p=3$, this can be attributed to the fact that removing the knots at $i=p+2$ and $i=\numSplines$ yields exactly the averaged knot vector when using $\numSplines_\mr{avg} = \numSplines_\mr{eq}-2$ in \eqref{eq:eqknots} and \eqref{eq:avgknots}.
Based on empirical investigations, we found that analogous boundary-knot removal strategies also stabilize the interpolation for higher polynomial degrees $p$.
The corresponding numbers of knots to remove, denoted by $\numSplines_\text{remove}$, are summarized in \Cref{tab:knotstoremove}.
\begin{table}
    \centering
    \caption{RECOMMENDED NUMBER OF KNOTS TO REMOVE FOR STABILIZED B-SPLINE INTERPOLATION.}
    \label{tab:knotstoremove}
    \begin{tabular}{c c}
    \toprule
    $p$ & Knots to remove \\
    \cmidrule(lr){1-1} \cmidrule(lr){2-2}
    1 & 0 \\
    2 & 0 \\
    3 & 1 \\
    4 & 1 \\
    5 & 2 \\
    6 & 2 \\
    \bottomrule
    \end{tabular}
\end{table}
For $p \neq 3$, the resulting knot vectors are no longer averaged knot vectors, but the removal still substantially reduces the growth of the stability constant, as shown in \Cref{fig:cond_stable}.
We refer to the resulting knot vectors as stabilized knot vectors.

The main advantage of this approach is that all knots of the stabilized knot vector are contained in the original equidistant knot vector.
Consequently, the mapping between the corresponding spline configurations can be constructed exactly via the knot-insertion procedure described in \Cref{ssec:bsplinetransfer}.
This allows the moments $\vec m_{\,d,\veg f_{\!n}}$ to be computed with respect to the original equidistant knot vector, thereby preserving exact and recursive interlevel transfers.
At the same time, the kernel interpolation is performed using the stabilized splines.

Specifically, stabilization transfer matrices $\mat S$ are used to map the moments from the equidistant splines to the stabilized splines.
These matrices are constructed using the same knot-insertion procedure as the interlevel transfer matrices $\mat C^\text{1D}$ in \Cref{ssec:bsplinetransfer}.
The stabilized interpolation weights are computed by
\begin{equation}
    \mat W^{\mr{B\text{-}spline}}_\text{stab} =
    \mat S^\T
    \left(\mat B_{\numSplines-\numSplines_\text{remove}, \text{stab}}^{\text{3D},p}\right)^{-\T}
    \mat G
    \left(\mat B_{\numSplines-\numSplines_\text{remove}, \text{stab}}^{\text{3D},p}\right)^{-1}
    \mat S,
\end{equation}
where $\mat B_{\numSplines-\numSplines_\text{remove},\text{stab}}^{\text{3D},p}$ denotes the stabilized sampling matrix.

The proposed stabilization strategy simultaneously satisfies all requirements of the multilevel scheme: the original equidistant knot vectors enable exact and recursively applicable interlevel transfers, the stabilized knot vectors provide numerically stable kernel interpolation, and the mapping between both spline configurations can itself be performed exactly via knot insertion without introducing additional errors.
The associated trade-off is that the effective number of splines available for interpolation is reduced by $2\numSplines_\text{remove}$.

\section{Factorizations of the Electric Field Integral Equation} \label{sec:factorization}
The proposed fast integral method relies on a factorization of the discretized \ac{EFIE} to enable efficient kernel interpolation and translation operations.
However, this factorization is not unique, since the derivatives in $\TPhi$ give rise to several possibilities for applying the interpolation to the \ac{EFIE} kernel.
The different factorizations affect the number and type of moments that must be computed, the storage requirements, and the computational cost of the translations.

In the following, we consider three different factorization strategies: a star-based factorization, a direct factorization, and a dyadic factorization.
The direct and dyadic formulations have been studied in~\cite{schobertLowFrequencySurfaceIntegral2012} previously.
In contrast, the star-based factorization is newly introduced in this work and extends the fast integral method by incorporating ideas from loop-star decompositions.
The resulting factorizations differ in their moment structure, storage requirements, and translation costs, as discussed in the subsequent sections.

\subsection{Star-Based Factorization}
\begin{figure}
    \centering
    \includegraphics{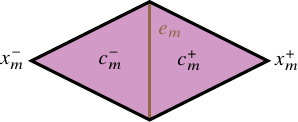}
    \caption{Geometric conventions for defining an \ac{RWG} function associated with the $m$th edge $e_m$, its adjacent cells $c_m^\pm$, and the corresponding vertices $x_m^\pm$.}
    \label{fig:plotrwg}
\end{figure}

For the star-based factorization, we propose to exploit the representation known from loop-star stabilization schemes of the discretized scalar potential operator \eqref{eq:matTPhi} as
\begin{equation}
    \matTPhi = \mat\Sigma \matV \mat\Sigma^{\mr T},
\end{equation}
where $\matV$ denotes the single-layer operator
\begin{equation}
    (\op V \phi)(\veg x) = \int_{\Gamma} g(\veg x, \veg y)\, \phi(\veg y)\, \mathrm dS(\veg y),
\end{equation}
discretized with piecewise constant basis functions
\begin{equation}
    \phi_m(\veg x) =
    \begin{cases}
    1/a_m, & \veg x \in c_m, \\
    0, & \text{otherwise},
    \end{cases}
\end{equation}
where $a_m$ denotes the area of the triangle $c_m$.
The sparse matrix $\mat\Sigma$ is the star matrix~\cite{vecchiLoopstarDecompositionBasis1999,andriulliLoopStarLoopTreeDecompositions2012}, given by
\begin{equation}
    \sbr{\mat \Sigma}_{mn}\! =
    \begin{cases}
        1\,, & \text{for } c_n = c_m^+ \\
        -1\,, & \text{for } c_n = c_m^- \\
        0, &\text{otherwise}\,,
    \end{cases}
\end{equation}
where the triangles $c_m^\pm$ sharing the edge $e_m$ follow the orientation of the \ac{RWG} basis functions depicted in \Cref{fig:plotrwg}.
Instead of factorizing $\matTPhi$ directly, we factorize $\mat V$.
The vector potential is factorized analogously to \eqref{eq:TAfactor}.
In the star-based factorization, both \ac{RWG} and piecewise-constant moments are required,
\begin{align}
    \int_{\supp \veg f_{\!m}} \veg f_{\!m}(\veg x)\, \ell_{i,\veg c_{\!\veg x}}\!(\veg x)\, \dd S(\veg x),
    &\qquad i = 1,\dots,\numSplines^3, \\
    \int_{c_m} \phi_{m}(\veg x)\, \ell_{i,\veg c_{\!\veg x}}\!(\veg x)\, \dd S(\veg x),
    &\qquad i = 1,\dots,\numSplines^3.
\end{align}

A drawback of the star-based factorization is that two sets of near-interactions must be stored.
A key advantage is that the null space of $\matTPhi$ is preserved up to machine precision despite the approximate nature of fast integral methods.
This property is particularly beneficial for preconditioning strategies, such as Calderón preconditioning~\cite{andriulliMultiplicativeCalderonPreconditioner2008}.

\subsection{Direct Factorization}
The direct factorization treats the vector and scalar potentials simultaneously and directly interpolates the kernel.
In this formulation, two sets of moments are required, involving the \ac{RWG} basis functions and their surface divergence,~\cite{schobertLowFrequencySurfaceIntegral2012}
\begin{align}
    \int_{\supp \veg f_{\!m}} \veg f_{\!m}(\veg x)\, \ell_{i,\veg c_{\!\veg x}}\!(\veg x)\, \dd S(\veg x),
    &\qquad i = 1,\dots,\numSplines^3, \\
    \int_{\supp \veg f_{\!m}} \Div_\Gamma \veg f_{\!m}(\veg x)\, \ell_{i,\veg c_{\!\veg x}}\!(\veg x)\, \dd S(\veg x),
    &\qquad i = 1,\dots,\numSplines^3.
\end{align}
As in the star-based factorization, this leads to four moment components per interpolation point: three vector components and one scalar divergence component.

\subsection{Dyadic Factorization}
The dyadic factorization is based on the dyadic formulation of the \ac{EFIE}~\cite{schobertLowFrequencySurfaceIntegral2012},
\begin{multline}
    \sbr{\mat T}_{mn} = \jm k \int_{\Gamma} \int_{\Gamma} \veg f_{\!m}(\veg x) \cdot \left(\veg I + \frac{1}{k^2} \grad_{\veg x} \grad_{\veg x} \!g(\veg x, \veg y) \right) \\ \veg f_{\!n}(\veg y)
    \, \dd S(\veg y) \dd S(\veg x),
\end{multline}
where $\veg I$ denotes the unit dyad.
In this formulation, only \ac{RWG} moments are required~\cite{schobertLowFrequencySurfaceIntegral2012},
\begin{equation}
    \int_{\supp \veg f_{\!m}}
    \veg f_{\!m}(\veg x)\, \ell_{i,\veg c_{\!\veg x}}\!(\veg x)\,
    \mathrm dS(\veg x),
    \qquad i = 1,\dots,\numSplines^3,
\end{equation}
which reduces storage requirements for the moments.
However, the translations become dyadic, with each translation involving 9 components corresponding to the dyad's entries, increasing computational complexity.

\section{Numerical Results} \label{sec:numerics}
\begin{figure}[tp]
    \centering
    \subfloat{\includegraphics[width=0.25\linewidth]{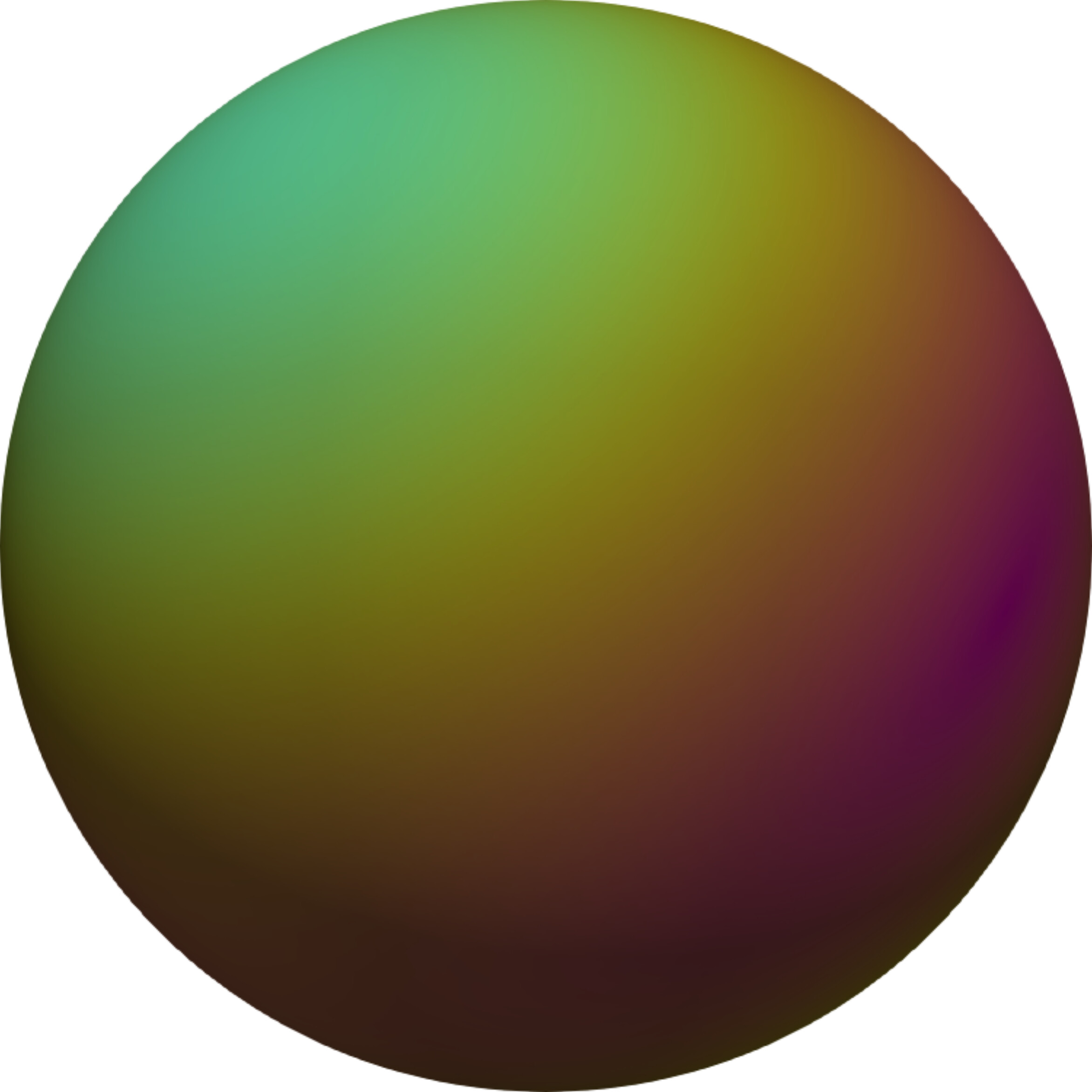}}
    \hspace{0.4cm}
    \subfloat{\includegraphics[width=0.25\linewidth]{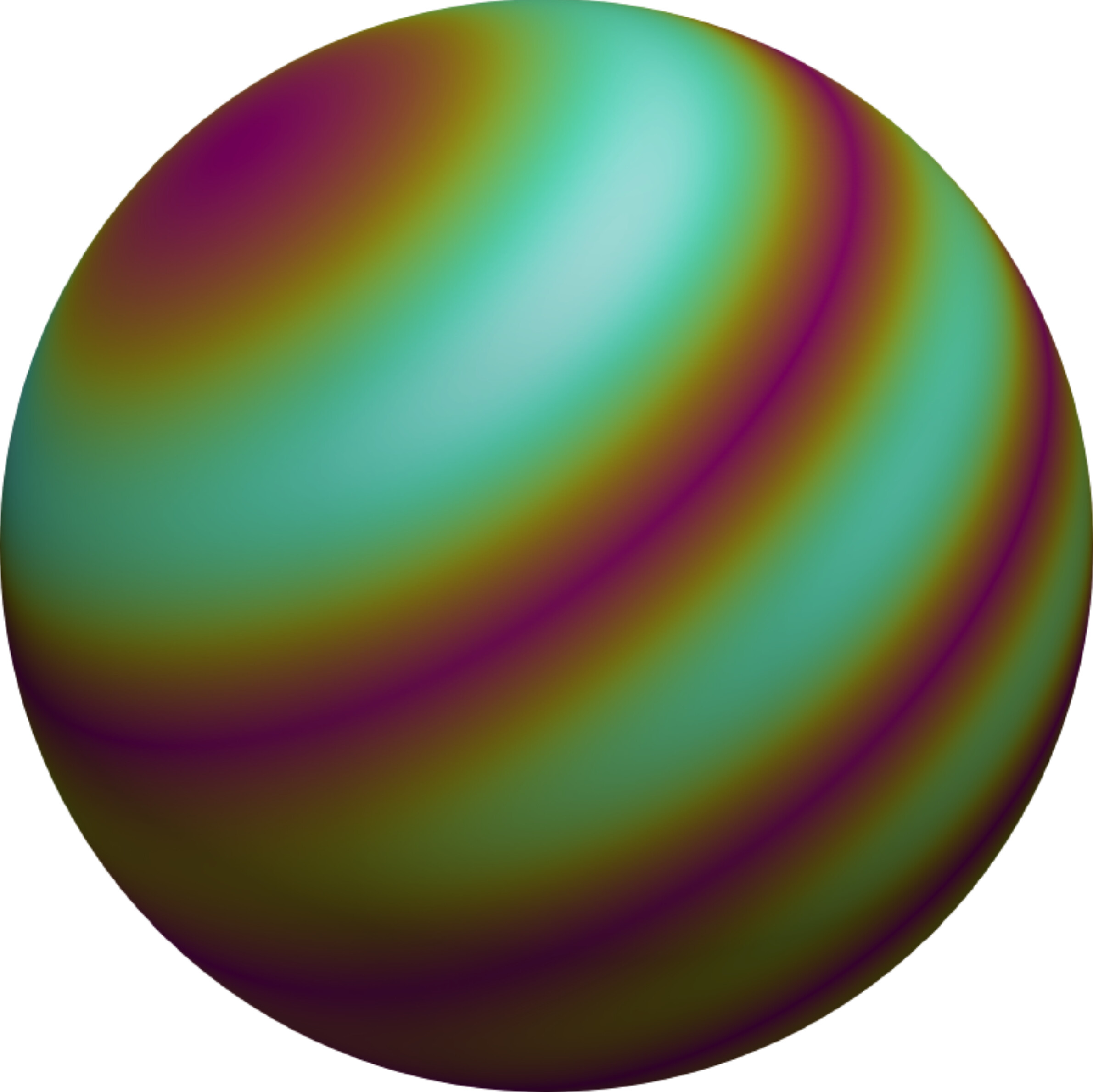}}
    \hspace{0.4cm}
    \subfloat{\includegraphics[width=0.25\linewidth]{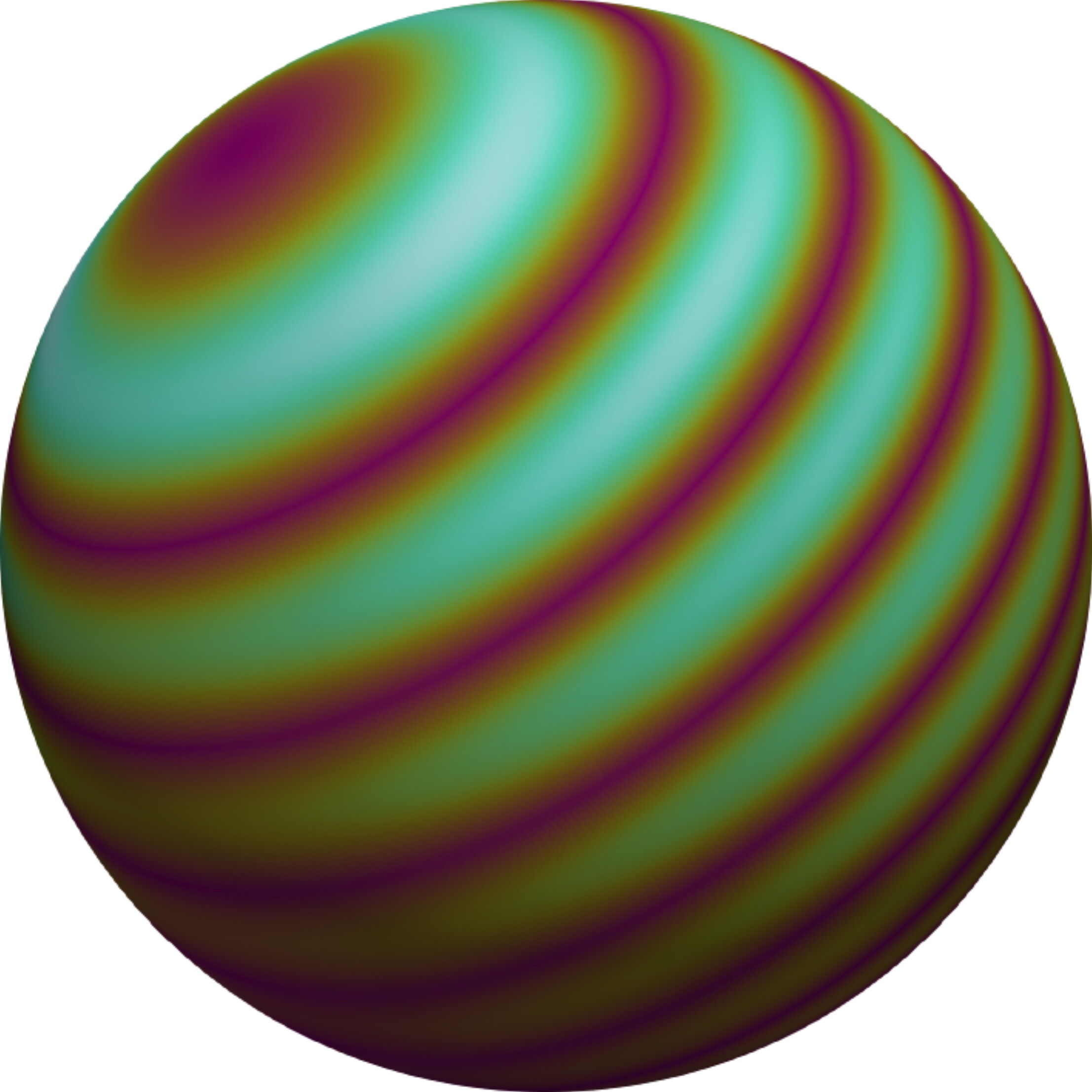}}\\
    \vspace{0.2cm}
    \subfloat{\includegraphics{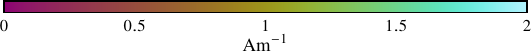}}
    \caption{Imaginary part of the surface current density $\veg j$ on the spheres excited by a plane wave. From left to right: diameters of $0.2\lambda$, $2\lambda$, and $4\lambda$.}
    \label{fig:plane_wave_spheres}
\end{figure}

\begin{figure*}[tp]
    \centering
    \subfloat{\includegraphics[width=0.3\linewidth]{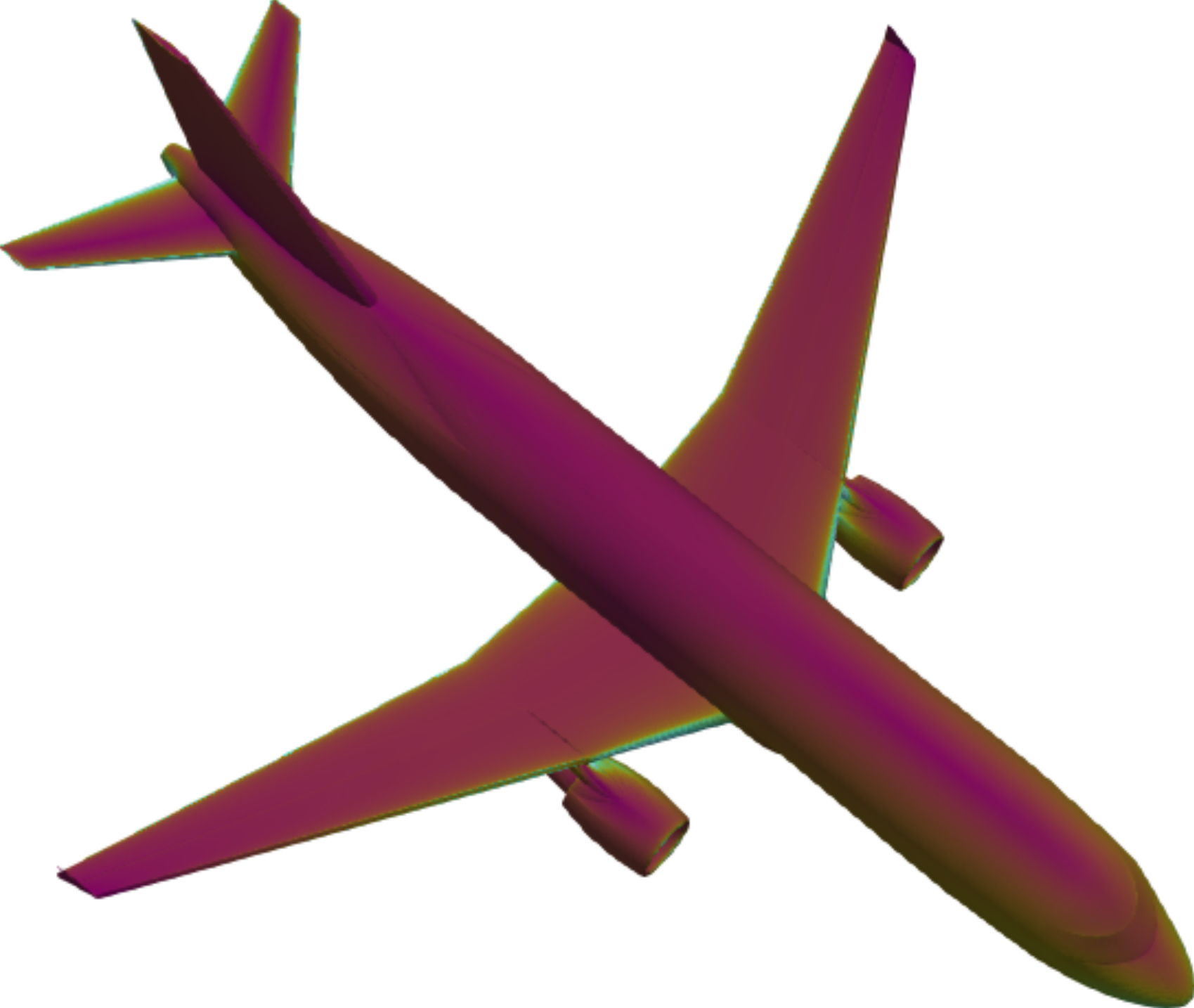}}
    \hspace{0.4cm}
    \subfloat{\includegraphics[width=0.3\linewidth]{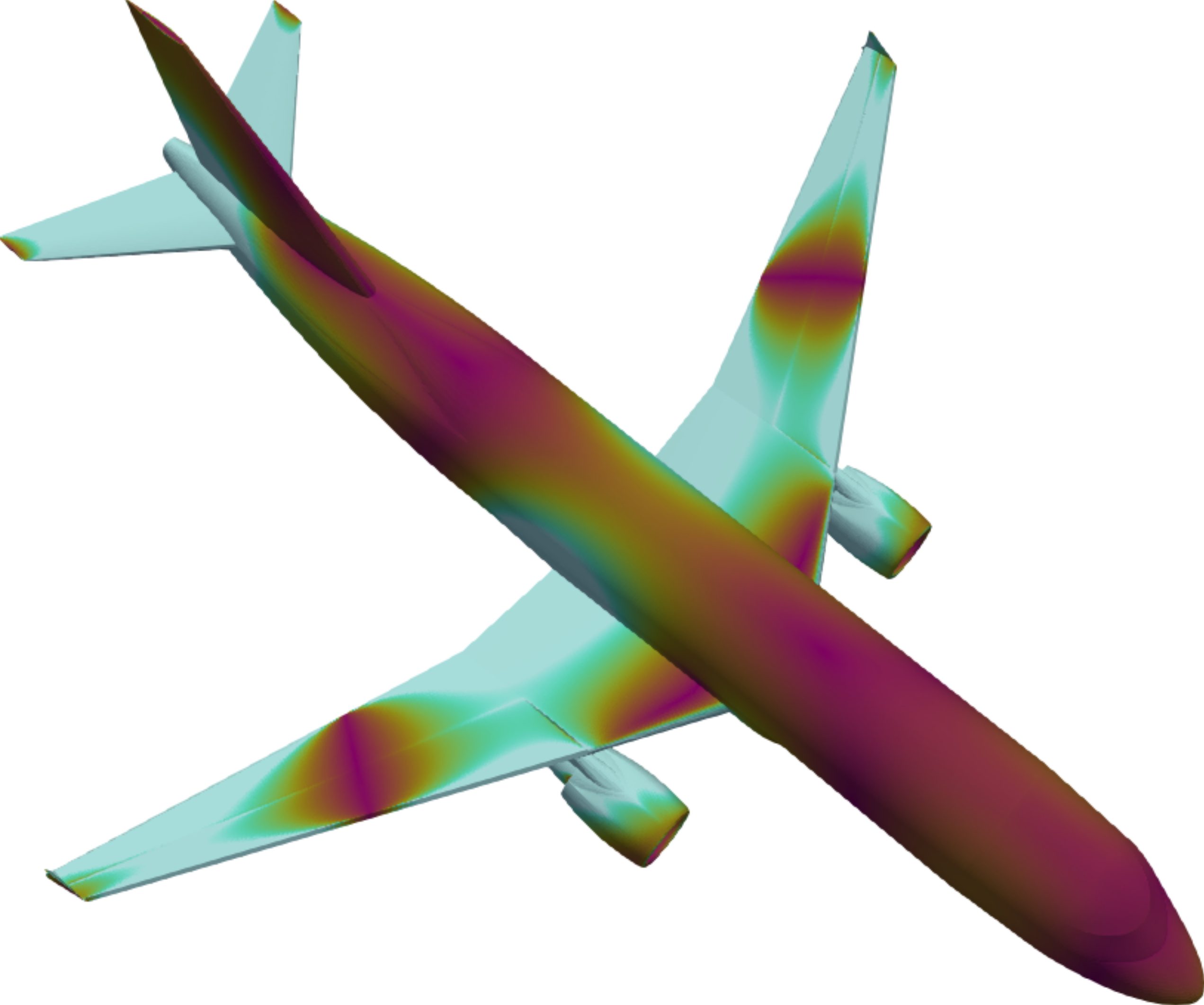}}
    \hspace{0.4cm}
    \subfloat{\includegraphics[width=0.3\linewidth]{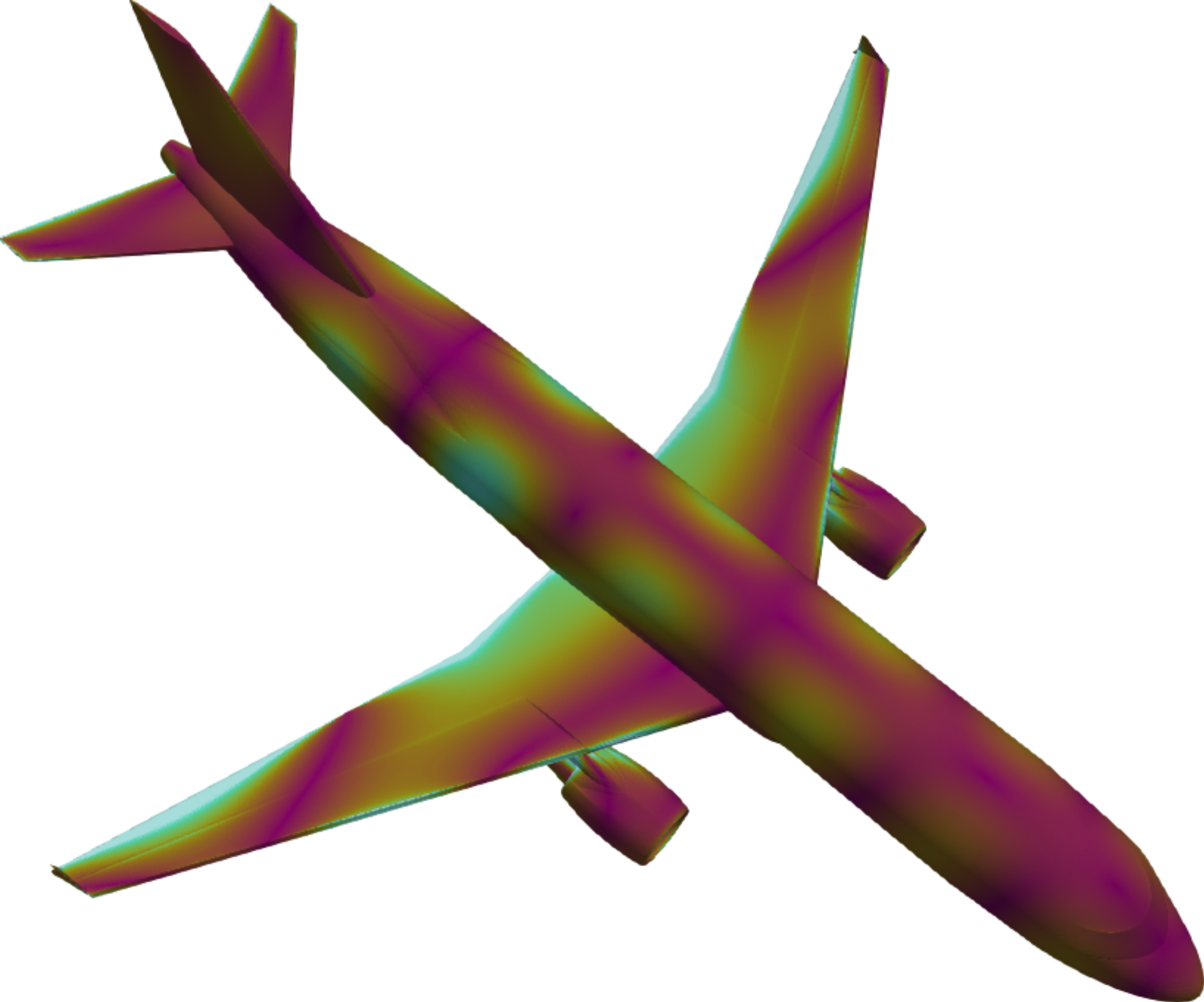}}\\
    \vspace{0.2cm}
    \subfloat{\includegraphics{figures/plot_sphere_plane_wave_colormap.pdf}}
    \caption{Imaginary part of the surface current density $\veg j$ on an aircraft model excited by a plane wave. From left to right: lengths of $0.2\lambda$, $2\lambda$, and $4\lambda$.}
    \label{fig:plane_wave_airbus}
\end{figure*}

To evaluate the proposed B-spline-based fast integral method with respect to numerical stability, accuracy, and computational performance, we consider both a canonical geometry (a sphere) and a realistic geometry (an aircraft model).
The solutions of the plane-wave excited scattering problems are shown in \Cref{fig:plane_wave_spheres,fig:plane_wave_airbus}, with discretizations comprising $\num{3168750}$ unknowns for the sphere and $\num{2028982}$ unknowns for the aircraft.

All simulations employ the multiplicative Calderón preconditioner~\cite{andriulliMultiplicativeCalderonPreconditioner2008} to reduce the number of iterations required by the iterative solver.

We first compare the accuracy of the Lagrange- and B-spline-based formulations for electrical sizes of $0.2\lambda$, $2\lambda$, and $4\lambda$.
Subsequently, the computational performance is analyzed in terms of \ac{MVP} time, memory consumption, and assembly time for the star-based factorization, focusing on a representative size of $0.2\lambda$.

A nonuniform octree is constructed by subdividing a box if it contains more than $\num{200}$ elements.
Refinement is terminated when a basis function would protrude more than $25\%$ of the box size.
For boxes with an edge length exceeding $0.5\lambda$, the number of interpolation polynomials is increased to maintain accuracy.

\subsection{Error Controllability} \label{sec:stabilityvalidation}
The stability of the fast integral methods (both the existing Lagrange-based compression and the proposed B-spline-based compression) is assessed using the relative spectral-norm error
\begin{equation}
    \aleph \left(\mat T \right) =\dfrac{\norm{\mat T^{\mathscr{H}^2,\mr{far}} - \mat T^\mr{far}}_2}{\norm{\mat T^\mr{far}}_2},
\end{equation}
which compares the far-interactions of the fully assembled matrix $\mat T$ and its $\mathscr{H}^2$-matrix approximation $\mat T^{\mathscr{H}^2}$.
The norm is evaluated using a power iteration method~\cite{misesPraktischeVerfahrenGleichungsauflosung1929}, and only far interactions are included in the error measure.

The error is analyzed for varying interpolation orders $p$ in $\mat T^{\mathscr{H}^2,\mr{far}}$ and, in the case of the B-spline interpolation, also in the number of splines $\numSplines$.
For Lagrange interpolation, we consider both strategies of handling protruding basis functions, as discussed in \Cref{sec:protruding}.

To obtain a more detailed assessment, the compression errors of $\matTA$ and $\matTPhi$ are evaluated separately, preventing potential error cancellation effects between the two system matrices.
A sphere and an aircraft model serve as representative test cases, discretized with $\num{48000}$ and $\num{62144}$ \ac{RWG} basis functions, respectively.

\subsubsection{Scalar Potential}
\begin{figure*}
    \centering

    \includegraphics{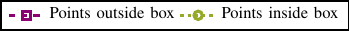}
    \vspace{2mm}

    \subfloat{\includegraphics{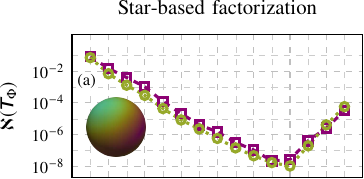}}
    \hfill
    \subfloat{\includegraphics{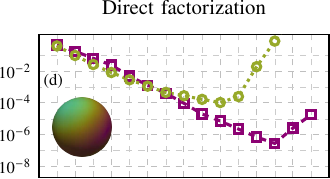}}
    \hfill
    \subfloat{\includegraphics{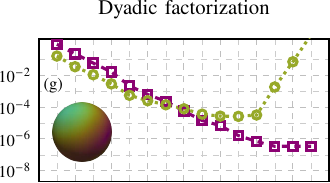}}\\[0.1cm]

    \subfloat{\includegraphics{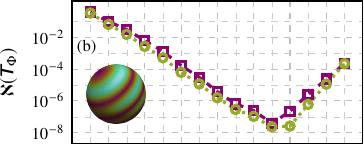}}
    \hfill
    \subfloat{\includegraphics{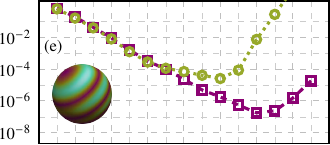}}
    \hfill
    \subfloat{\includegraphics{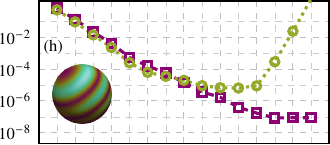}}\\[0.1cm]

    \subfloat{\includegraphics{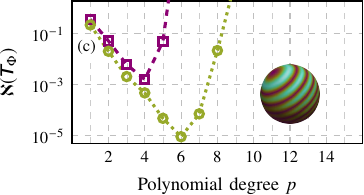}}
    \hfill
    \subfloat{\includegraphics{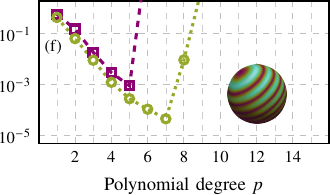}}
    \hfill
    \subfloat{\includegraphics{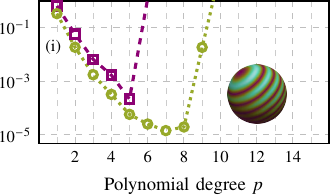}}

    \caption{
    Relative spectral-norm error between the Lagrange-based fast approximation $\matTPhi^{\mathscr{H}^2}$ and the fully assembled matrix $\matTPhi$ for a sphere discretized with $\num{48000}$ \ac{RWG} functions.
    Columns correspond to the star-based ((a)--(c)), direct ((d)--(f)), and dyadic ((g)--(i)) factorizations, while rows correspond to spheres with diameters $0.2\lambda$, $2\lambda$, and $4\lambda$.
}
    \label{fig:lagrange_sphere_scalarpotential_all}
\end{figure*}

We first consider, as a canonical example, the case of a sphere, for which we analyze the star-based, direct, and dyadic factorizations introduced in \Cref{sec:factorization}.
The results for the Lagrange-based compression are shown in \Cref{fig:lagrange_sphere_scalarpotential_all}.

As expected, the Runge phenomenon leads to pronounced instabilities.
These instabilities occur at lower interpolation orders as the electrical size increases, since larger boxes require higher polynomial degrees to maintain accuracy.
For the smaller spheres ($0.2\lambda$ and $2\lambda$), placing interpolation points outside the box improves the accuracy of the direct and dyadic factorizations.
However, for the largest sphere ($4\lambda$), this strategy causes an earlier onset of instability.

\begin{figure*}
    \centering

    \includegraphics{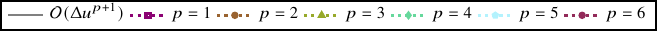}
    \vspace{2mm}

    \subfloat{\includegraphics{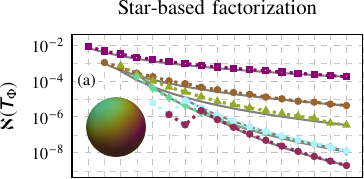}}
    \hfill
    \subfloat{\includegraphics{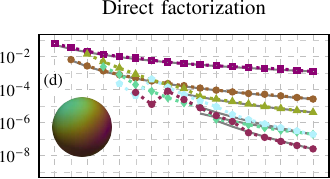}}
    \hfill
    \subfloat{\includegraphics{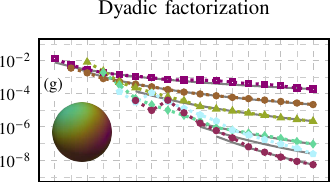}}\\[0.1cm]

    \subfloat{\includegraphics{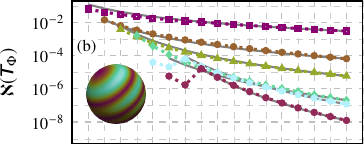}}
    \hfill
    \subfloat{\includegraphics{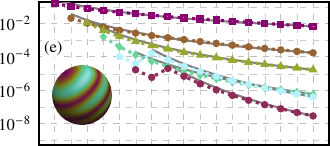}}
    \hfill
    \subfloat{\includegraphics{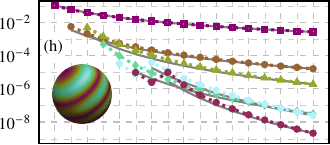}}\\[0.1cm]

    \subfloat{\includegraphics{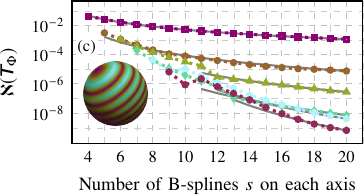}}
    \hfill
    \subfloat{\includegraphics{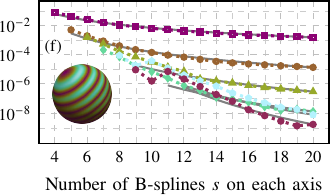}}
    \hfill
    \subfloat{\includegraphics{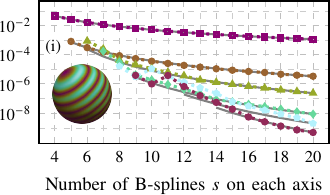}}

    \caption{
    Relative spectral-norm error between the B-spline-based fast approximation $\matTPhi^{\mathscr{H}^2}$ and the fully assembled matrix $\matTPhi$ for a sphere discretized with $\num{48000}$ \ac{RWG} functions.
    Columns correspond to the star-based ((a)--(c)), direct ((d)--(f)), and dyadic ((g)--(i)) factorizations, while rows correspond to spheres with diameters $0.2\lambda$, $2\lambda$, and $4\lambda$.
}    \label{fig:sphere_scalarpotential_all}
\end{figure*}

In contrast, \Cref{fig:sphere_scalarpotential_all} demonstrates that the B-spline-based kernel interpolation does not suffer from the Runge phenomenon.
The results exhibit convergence behavior consistent with the asymptotic error bound~\eqref{eq:errorapprox} derived in \Cref{sec:errorestimates} as $\numSplines$ and $p$ increase.
Compared to the Lagrange-based method, lower errors are obtained for a sufficiently high combination of $\numSplines$ and $p$.
The jumps in the error observed for $p=5$ and $p=6$ are caused by the stabilization strategy described in \Cref{ssec:stab}, which removes splines from the interpolation to improve numerical stability.

\begin{figure}
    \centering
    \makebox[\columnwidth][r]{%
        \subfloat[]{
            \includegraphics{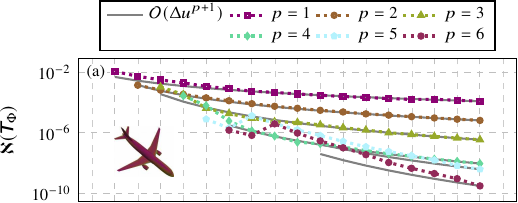}
            \label{fig:airbus_scalarpotential_0.1}
        }
    }\\[-0.4cm]
    \makebox[\columnwidth][r]{%
        \subfloat[]{
            \includegraphics{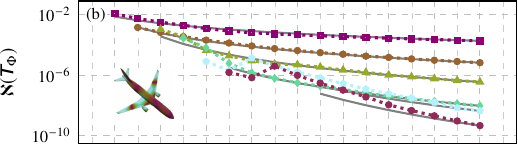}
            \label{fig:airbus_scalarpotential_1}
        }
    }\\[-0.4cm]
    \makebox[\columnwidth][r]{%
        \subfloat[]{
            \includegraphics{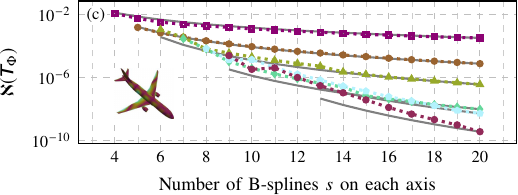}
            \label{fig:airbus_scalarpotential_2}
        }
    }\\[-0.4cm]
    \caption{
    Relative spectral-norm error between the B-spline-based fast approximation $\matTPhi^{\mathscr{H}^2}$ and the fully assembled matrix $\matTPhi$ for an aircraft model discretized with $\num{62144}$ \ac{RWG} functions using the star-based factorization.
    Rows correspond to models with lengths $0.2\lambda$, $2\lambda$, and $4\lambda$.
    }
    \label{fig:airbus_scalarpotential}
\end{figure}

Among the three factorizations, the star-based factorization exhibits the most robust error behavior.
Thus, we limit ourselves to this case for study, as an example of a realistic geometry, the aircraft model.
As shown in \Cref{fig:airbus_scalarpotential}, the error also decreases with increasing $\numSplines$ and $p$ in this case.
Overall, increasing the polynomial degree $p$ reduces the interpolation error more effectively than increasing the number of splines $\numSplines$, in agreement with \eqref{eq:errorapprox}.
The number of splines $\numSplines$ provides an additional tuning parameter that can further reduce the error, particularly in the high-accuracy regime.

\subsubsection{Vector Potential}

\begin{figure}
    \centering
    \makebox[\columnwidth][r]{%
        \subfloat[]{
            \includegraphics{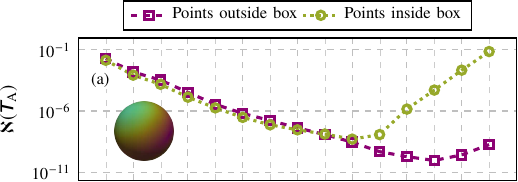}
            \label{fig:lagrange_sphere_vectorpotential_0.1}
        }
    }\\[-0.4cm]
    \makebox[\columnwidth][r]{%
        \subfloat[]{
            \includegraphics{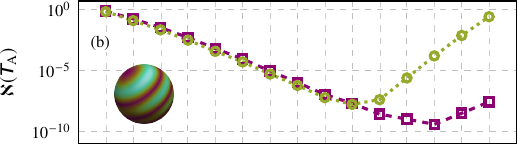}
            \label{fig:lagrange_sphere_vectorpotential_1}
        }
    }\\[-0.4cm]
    \makebox[\columnwidth][r]{%
        \subfloat[]{
            \includegraphics{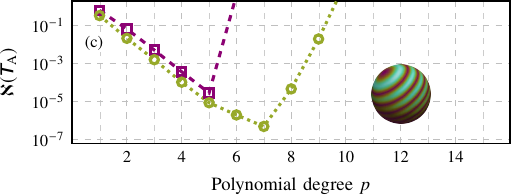}
            \label{fig:lagrange_sphere_vectorpotential_2}
        }
    }\\[-0.4cm]
    \caption{
    Relative spectral-norm error between the Lagrange-based fast approximation $\matTA^{\mathscr{H}^2}$ and the fully assembled matrix $\matTA$ for a sphere discretized with \num{48000} \ac{RWG} functions.
    Rows correspond to spheres with diameters $0.2\lambda$, $2\lambda$, and $4\lambda$.
    }
    \label{fig:lagrange_sphere_vectorpotential}
\end{figure}

\begin{figure}
    \centering
    \makebox[\columnwidth][r]{%
        \subfloat[]{
            \includegraphics{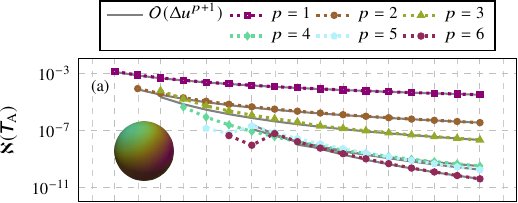}
            \label{fig:sphere_vectorpotential_0.1}
        }
    }\\[-0.4cm]
    \makebox[\columnwidth][r]{%
        \subfloat[]{
            \includegraphics{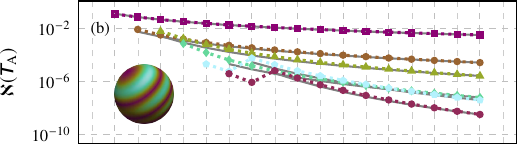}
            \label{fig:sphere_vectorpotential_1}
        }
    }\\[-0.4cm]
    \makebox[\columnwidth][r]{%
        \subfloat[]{
            \includegraphics{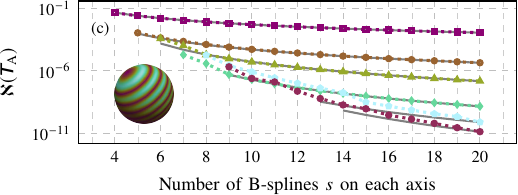}
            \label{fig:sphere_vectorpotential_2}
        }
    }\\[-0.4cm]
    \caption{
    Relative spectral-norm error between the B-spline-based fast approximation $\matTA^{\mathscr{H}^2}$ and the fully assembled matrix $\matTA$ for a sphere discretized with \num{48000} \ac{RWG} functions.
    Rows correspond to spheres with diameters $0.2\lambda$, $2\lambda$, and $4\lambda$.
    }
    \label{fig:sphere_vectorpotential}
\end{figure}

\begin{figure}
    \centering
    \makebox[\columnwidth][r]{%
        \subfloat[]{
            \includegraphics{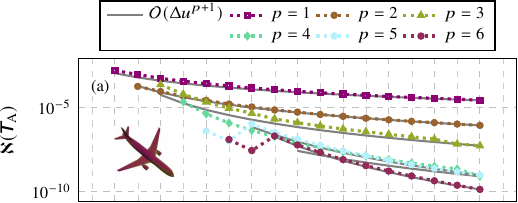}
            \label{fig:airbus_vectorpotential_0.1}
        }
    }\\[-0.4cm]
    \makebox[\columnwidth][r]{%
        \subfloat[]{
            \includegraphics{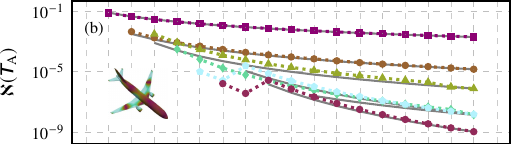}
            \label{fig:airbus_vectorpotential_1}
        }
    }\\[-0.4cm]
    \makebox[\columnwidth][r]{%
        \subfloat[]{
            \includegraphics{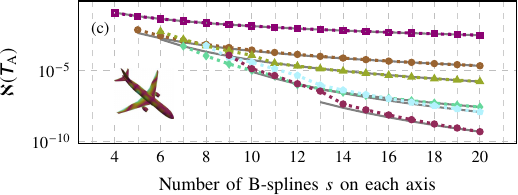}
            \label{fig:airbus_vectorpotential_2}
        }
    }\\[-0.4cm]
  \caption{
    Relative spectral-norm error between the B-spline-based fast approximation $\matTA^{\mathscr{H}^2}$ and the fully assembled matrix $\matTA$ for an aircraft model discretized with $\num{62144}$ \ac{RWG} functions.
    Rows correspond to models with lengths $0.2\lambda$, $2\lambda$, and $4\lambda$.
    }
    \label{fig:airbus_vectorpotential}
\end{figure}

Results for the vector potential are shown for the sphere in \Cref{fig:lagrange_sphere_vectorpotential,fig:sphere_vectorpotential} and for the aircraft model in \Cref{fig:airbus_vectorpotential}.

Again, we first study the case of a sphere.
As for the scalar potential, \Cref{fig:lagrange_sphere_vectorpotential} shows that the Lagrange-based method breaks down as the polynomial degree increases due to the Runge phenomenon.
However, overall lower errors are achieved compared with compressing the scalar potential.
In contrast, \Cref{fig:sphere_vectorpotential} displays that the B-spline-based method remains stable for all tested electrical sizes.
As in the scalar-potential case, jumps in the error are observed for $p=5$ and $p=6$ due to the stabilization strategy described in \Cref{ssec:stab}.
\Cref{fig:airbus_vectorpotential} confirms that the compression remains stable in the case of a realistic model.

\subsection{Performance Analysis}
All computations are performed on a sysGen Server G493-ZB1-AAP1 equipped with two dual AMD EPYC TURIN 9755, each with 128 cores and a clock frequency of $\SIrange{2.70}{4.10}{GHz}$, and $\SI{6}{TB}$ DDR5 RAM.
The simulations use 16 threads in single-precision arithmetic, whereas all error estimates are computed in double-precision arithmetic to accurately assess the interpolation error.
We evaluate the runtime of the \ac{MVP}, memory consumption, assembly time, and relative error of the B-spline-based fast integral method for both the sphere and aircraft geometries as a function of $N$ in order to assess the predicted $\mc O(N)$ asymptotic complexity.
Unless stated otherwise, the interpolation parameters are chosen as $\numSplines =p+4$.
This choice yields a knot spacing of $\Updelta u = 1/4$, such that a single buffer spline is sufficient to cover the maximum protrusion of $25\%$ of the box size while still leaving half a box distance between interpolation domains.
Empirically, this was found to provide a good trade-off between accuracy and robustness.
All experiments employ the proposed star-based factorization.

To estimate the compression error introduced by the B-spline-based fast integral method, we use a second B-spline-based fast integral method $\mat T_\mr{ref}^{\mathscr{H}^2}$ with $p=6$ and $\numSplines=14$ as a reference solution and compute the relative error
\begin{equation}
    \aleph_\mr{self} \left(\mat T^{\mathscr{H}^2} \right) = \dfrac{\norm{\mat T^{\mathscr{H}^2,\mr{far}}-\mat T_\mr{ref}^{\mathscr{H}^2,\mr{far}}}_2}{\norm{\mat T_\mr{ref}^{\mathscr{H}^2,\mr{far}}}_2
    }\,,
\end{equation}
where only the far interactions are considered.

We first consider, as a canonical example, a sphere with diameter $0.2\lambda$.
Since no qualitative differences were observed for other diameters, these results are omitted.
The number of unknowns $N$ ranges from $\num{12000}$ to $\num{3072000}$.

\begin{figure}
    \centering

    \makebox[\columnwidth][r]{%
        \subfloat[]{
            \includegraphics{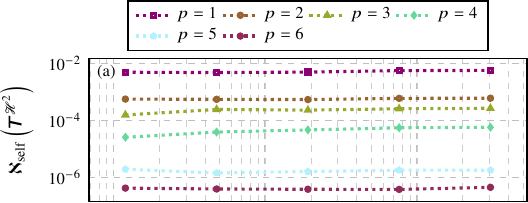}
            \label{fig:sphere_err}
        }
    }\\[-0.4cm]
    \makebox[\columnwidth][r]{%
        \subfloat[]{
            \includegraphics{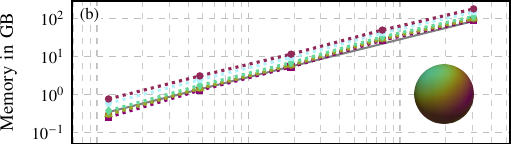}
            \label{fig:sphere_storage}
        }
    }\\[-0.4cm]
    \makebox[\columnwidth][r]{%
        \subfloat[]{
            \includegraphics{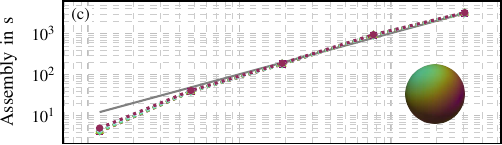}
            \label{fig:sphere_time}
        }
    }\\[-0.4cm]
    \makebox[\columnwidth][r]{%
        \subfloat[]{
            \includegraphics{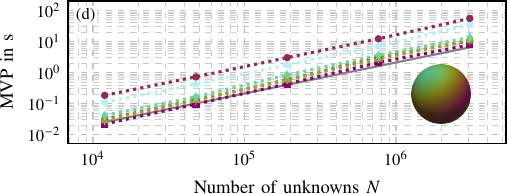}
            \label{fig:sphere_mv}
        }
    }\\[-0.4cm]
    \caption{
    Sphere: performance analysis of the B-splined based compression with star-based factorization for a sphere with diameter $0.2\lambda$.
    (a) Relative error.
    (b) Storage.
    (c) Assembly time.
    (d) Time of single \ac{MVP}.}

    \label{fig:timings_storage_assembly_sphere}
\end{figure}

\Cref{fig:sphere_err} confirms that the compression of $\mat T$ remains well controlled as $N$ increases.
Increasing the polynomial degree $p$ consistently reduces the error, in agreement with the error estimate in \eqref{eq:errorapprox}.
However, the reduction in error from $p=2$ to $p=3$ is less pronounced because the stabilization strategy described in \Cref{ssec:stab} removes splines from the interpolation.

\Crefrange{fig:sphere_storage}{fig:sphere_mv} confirm the $\mc O(N)$ scaling for storage, assembly time, and \ac{MVP} cost.
For the discretization with $N=\num{3072000}$, \Cref{fig:sphere_storage} shows that the memory consumption increases from \SI{91}{GB} for $p=1$ to \SI{182}{GB} for $p=6$.
This increase is expected, since the number of splines is chosen as $\numSplines=p+4$.
In contrast, the assembly time is barely affected by $p$ as shown in \Cref{fig:sphere_time}, where the time ranges from $\SI{51}{min}$ to $\SI{53}{min}$ for $N=\num{3072000}$.
In fact, the assembly is dominated by the computation of the near interactions, which account for at least $79\%$ of the total assembly time.
The \ac{MVP} cost exhibit a stronger dependence on $p$, ranging from $\SI{7.6}{s}$ to $\SI{54}{s}$ for $N=\num{3072000}$.

\begin{figure}
    \centering

    \makebox[\columnwidth][r]{%
        \subfloat[]{
            \includegraphics{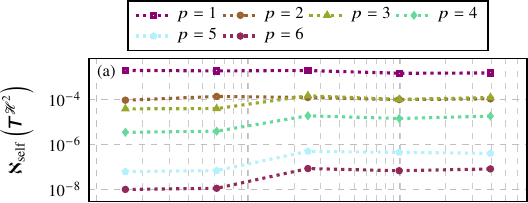}
            \label{fig:aircraft_err}
        }
    }\\[-0.4cm]
    \makebox[\columnwidth][r]{%
        \subfloat[]{
            \includegraphics{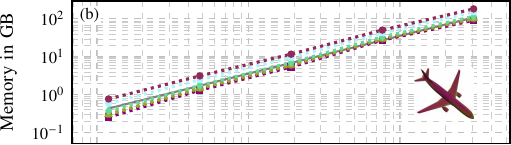}
            \label{fig:aircraft_storage}
        }
    }\\[-0.4cm]
    \makebox[\columnwidth][r]{%
        \subfloat[]{
            \includegraphics{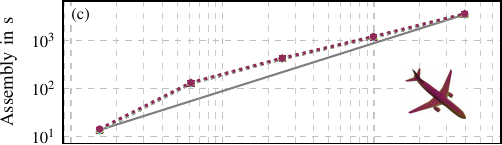}
            \label{fig:aircraft_time}
        }
    }\\[-0.4cm]
    \makebox[\columnwidth][r]{%
        \subfloat[]{
            \includegraphics{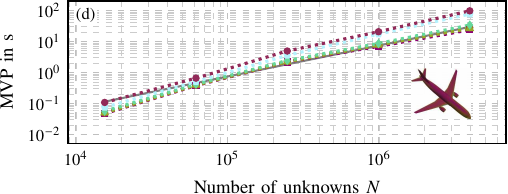}
            \label{fig:aircraft_mv}
        }
    }\\[-0.4cm]
    \caption{
    Aircraft: performance analysis of the B-splined based compression with star-based factorization for an aircraft with diameter $0.2\lambda$.
    (a) Relative error.
    (b) Storage.
    (c) Assembly time.
    (d) Time of single \ac{MVP}.}

    \label{fig:timings_airbus}
\end{figure}

As a realistic multiscale example, we consider an aircraft geometry, where the number of unknowns ranges from $\num{15532}$ to $\num{3978112}$.
\Cref{fig:aircraft_err} indicates that the compression error remains well controlled as $N$ increases.
For polynomial degrees $p=3$ to $p=6$, the error increases moderately from $N=\num{248608}$ onward before saturating.

\Crefrange{fig:aircraft_storage}{fig:aircraft_mv} confirm the expected $\mc O(N)$ scaling for storage, assembly time, and \ac{MVP} cost.
For $N=\num{3978112}$, the memory consumption ranges from \SI{347}{GB} to \SI{439}{GB}, the assembly time from \SI{57}{min} to \SI{60}{min}, and the cost of a \ac{MVP} from \SI{23}{s} to \SI{95}{s} as $p$ increases from $1$ to $6$, demonstrating that the observed scaling behavior also carries over to realistic multiscale geometries.

\section{Conclusion}
Fast integral methods based on global Lagrange kernel interpolation are ultimately limited by the Runge phenomenon.
Neither extrapolation nor scaling of interpolation points eliminates this instability; both strategies merely shift the regime in which it appears.
The proposed stabilized B-spline interpolation avoids this breakdown while retaining exact interlevel transfers.
As a result, the method achieves substantially lower compression errors and remains stable even at high interpolation orders.
The newly proposed star-based factorization preserves the null space of the scalar-potential operator to machine precision, whereas the direct and dyadic factorizations offer alternative trade-offs between memory consumption and \ac{MVP} cost.
For electrically larger interactions, it can be combined with high-frequency techniques such as the \ac{MLFMM}, enabling broadband formulations that retain the advantages of stabilized B-spline interpolation across a wide range of electrical sizes~\cite{jukicBroadbandStabilizedMultilevel2025a}.
Overall, the results demonstrate that stabilized B-spline interpolation provides a practical high-accuracy alternative to Lagrange interpolation for fast integral-equation solvers.

\section*{Acknowledgment}
Boundary-element simulations were performed with the open-source Julia library BEAST.jl~\cite{BEAST}. B-spline computations were performed with the open-source Julia library NURBS.jl~\cite{nurbs2023}.
\ifCLASSOPTIONcaptionsoff
\newpage
\fi



\printbibliography

\end{document}